\documentclass[11pt]{article}
\usepackage[pctex32]{graphics}
\usepackage{amssymb}
\usepackage{latexsym}
\usepackage{epsfig}
\usepackage{floatflt}
\usepackage{multirow}
\usepackage{subfigure}
\usepackage{wrapfig}
\usepackage{mathrsfs}
\usepackage[left=1in,top=1in,bottom =1.4in,right=1in,nohead]{geometry}
\usepackage{pstcol}

\definecolor{Blue}{rgb}{0.,0.,0.6}
\definecolor{Green}{rgb}{0.,1.,0.0}
\definecolor{Red}{rgb}{1,0.,0.}
\definecolor{Purple}{rgb}{0.5,0.,0.5}

\newcommand{\R}{{\mathbb R}}

\newcommand{\mA}{{\mathsf A}}

\newcommand{\mD}{{\mathsf D}}
\newcommand{\mL}{{\mathsf L}}
\newcommand{\mI}{{\mathsf I}}
\newcommand{\mG}{{\mathsf G}}
\newcommand{\mF}{{\mathsf F}}
\newcommand{\mU}{{\mathsf U}}
\newcommand{\mV}{{\mathsf V}}
\newcommand{\mP}{{\mathsf P}}
\newcommand{\mW}{{\mathsf W}}
\newcommand{\mSigma}{{\mathsf \Sigma}}
\newcommand{\mS}{{\mathsf S}}
\newcommand{\mT}{{\mathsf T}}
\newcommand{\mM}{{\mathsf M}}
\newcommand{\mB}{{\mathsf B}}
\newcommand{\mR}{{\mathsf R}}
\newcommand{\mQ}{{\mathsf Q}}
\newcommand{\mO}{{\mathsf O}}
\newtheorem{definition}{Definition}[section]
\newtheorem{theorem}[definition]{Theorem}

\newtheorem{problem}[definition]{Problem}

\begin{document}

\title{Adaptive anisotropic Bayesian meshing for inverse problems}
\author{A Bocchinfuso\footnote{Current address: Oak Ridge National Labotatory} \and D Calvetti \and E Somersalo}
\date{Department of Mathematics, Applied Mathematics and Statistics \\
Case Western Reserve University}

\maketitle
\begin{abstract}
We consider inverse problems estimating distributed parameters from indirect noisy observations through discretization of continuum models described by partial differential or integral equations. It is well understood that the errors arising from the discretization can be detrimental for ill-posed inverse problems, as discretization error behaves as correlated noise. While this problem can be avoided with a discretization fine enough to suppress the modeling error level below that of the exogenous noise that is addressed, e.g., by regularization, the computational resources needed to deal with the additional degrees of freedom may require high performance computing environment. Following an earlier idea, we advocate the notion that the discretization is one of the unknowns of the inverse problem, and is updated iteratively together with the solution. In this approach, the discretization, defined in terms of an underlying metric, is refined selectively only where the representation power of the current mesh is insufficient. In this paper we allow the metrics and meshes to be anisotropic, and we show that this leads to significant reduction of memory allocation and computing time.
\end{abstract}

\section{Introduction}

The estimation of distributed parameters from indirect noisy observations is an inverse problem arising in a wide range of applications, including geophysics, non-destructive material evaluation and medical imaging. The distributed parameters may contain information about the underlying medium such as density, elastic parameters, electric conductivity or optical absorption and scattering coefficients. Typically these parameters are related to the observations through continuum models governed by partial differential equations or integral equations whose numerical solutions require discretization via finite differences, finite element, or finite volumes. To ensure that the discretized problem is close to the continuum one, the mesh must be fine enough for the discretization error not to hamper the approximation. The inherent ill-posedness of inverse problems means that small errors in data may be significantly amplified in the computed solutions, thus proper handling of the discretization error is even more critical in this context than in simulations or model predictions. Since numerical approximation errors arising from coarse discretization affect all output parameters simultaneously, they manifest themselves as if highly correlated noise had been added to the observations \cite{KSbook,KSmoderror,CDSS} and without a careful analysis of the correlation structure and suitable compensation, they may be detrimental for solving the inverse problem. Because the correlation of the discretization error is typically very structured, it is not well described in terms of the exogenous observation noise, and the higher the quality of the measured data is, the more noticeable the artifacts due to discretization in the computed solution may be.

There are two different ways to mitigate the effect of discretization errors. The first and most obvious approach is to refine the mesh until the error due to the discretization is smaller than the observation error. The drawback of this approach is that refining the discretization may increase the computational complexity of the forward model to the point where the numerical treatment of the inverse problem becomes impractical or impossible. An alternative is to model the approximation error in the Bayesian statistical framework by estimating its probability distribution and incorporating it in the likelihood model that acknowledges the imprecision of the forward model used for the inversion. The latter approach, which has been shown to consistently produce high quality reconstructions without excessive computational complexity \cite{Arridge,KSbook,KSmoderror,CDSS,Nissinen} once the density of the modeling error has been estimated, is not without a price. In fact, since the modeling error is a function of the unknown distributed parameter to be evaluated, estimating its distribution typically requires off-line nmerical forward simulations that may be very costly. Philosophically, the  Bayesian modeling error approach can be seen as a training algorithm similar to machine learning algorithms.

In this article, we advocate a third approach that reduces the discretization error by selectively refining the mesh. We assume that the distributed parameter is piecewise constant or nearly constant function, so that its generalized gradient admits a sparse or compressible approximation in terms of appropriately chosen basis functions. The inverse problem is solved by an iterative scheme in which the support of the sparse coefficient vector of the approximation of the gradient is learned while estimating the solution, and the computational domain is iteratively remeshed. The meshing is a function of an underlying metric that is updated at each iteration concomitantly with the approximate solution. The computational engine used for estimating the sparse representation of the gradient is an iterative numerical scheme based on hierarchical Bayesian hypermodels, a hybrid version of the iterative alternating scheme (IAS) studied extensively in the literature \cite{CPS,CPSS,CPPSV,CSS,Glaubitz1,Glaubitz2,PCS,Pursiainen}.

This article elaborates and extends the work in \cite{CCPS}, with a special attention to the computational aspects. One of the novelties of this work is that the approximation of  the generalized gradient of the unknown using Whitney elements \cite{Bossavit,Bossavit2} bridges the sparse representation and a geometric interpretation of it. Using the Whitney representation, we build a Cl\'{e}ment-type interpolant \cite{Brenner,Clement}, a continuous approximation of the generalized gradient that allows us to define the mesh-generating metric. This new formalism removes the need to introduce intermediate meshing interpolating the hyperparameter associated to the edges of the mesh that made the algorithm in \cite{CCPS} slightly cumbersome, and required some hand-tuning of the mesh near the domain boundary. With a continuous interpolant of the gradient at our disposal, we then define the  metric  in terms of the gradient. Another novelty is that now we can allow the metric and the corresponding mesh to be anisotropic, leading to a more efficient representation of the unknown, requiring less elements to approximate discontinuities. Moreover, instead of seeking a way to approximate the modeling error, the more efficient adaptive anisotropic meshing that we propose here is able to selectively refine the discretization so that the modeling error falls below the level of the observation noise. Finally, the Bayesian updating algorithm introduced in this article uses a hybrid version of the IAS algorithm, further improving the sparsity of the gradient representation, thus reducing the overall computing cost.

The computed examples presented in this article are similar to those considered in \cite{CCPS}, i.e., the X-ray fanbeam tomography problem with few illumination angles, and the inverse source problem for Darcy flow, thus facilitating the comparison of the novel algorithm  with the  adaptive isotropic Bayesian meshing.  Extending the Bayesian adaptive anisotropic meshing presented here to non-linear models, while in principle natural, will entail a further analysis of the solver, and will be is the topic of a separate contribution.

\section{Bridging the continuum and discrete formulations}

Consider a continuous observation model: Let $\Omega\subset\R^d$, $d=2,3$, be a bounded domain with a Lipschitz boundary, and let $u$ be a distributed parameter function defined in $\Omega$.  Let $F^*$ denote a function that models the dependency of a finite dimensional observable quantity $b$ on $u$, and consider an observation model with additive noise, 
\begin{equation}\label{cont model}
 b = F^*(u) + \varepsilon, \quad b\in\R^m, \quad u:\Omega\to\R.
\end{equation}
The inverse problem is to estimate $u$ based on observation of $b$. The forward model $F^*$  is typically  based on  PDE/integral equation. For the sake of definiteness, we assume that the function $u$ is a piecewise constant, or more generally, piecewise smooth function with small gradient outside the discontinuities, whose singular support we want to estimate numerically. 

Typically, when solving the inverse problem, the unknown $u$ as well as the forward model $F^*$ need to be discretized. 
 It is well understood \cite{KSmoderror} that the discrete approximation differs from the presumably accurate theoretical continuous model by a {\em model discrepancy} representing the discretization error. Assuming that the data are well described by the continuous model, the model discrepancy acts as an additional noise term. If the data are of high quality in the sense that the exogenous noise $\varepsilon$ is of low level, the model discrepancy may be the dominant part of the noise, and as inverse problems are ill-posed, ignoring it may lead to significant loss of accuracy in the estimate. Modeling the discrepancy, on the other hand, is a non-trivial task for several reasons:  First, the modeling error may have a non-zero mean and a complex correlation structure, so inflating artificially the noise level in the inverse solver may not have the desired effect. In fact, the model discrepancy depends on the discretization of the continuous model and on the unknown parameter $u$ itself, complicating the process of estimating it.

To address the modeling error question in a more rigorous framework, consider a discretization scheme of the forward model. For simplicity, we limit ourselves to the case $d=2$. Assume that the boundary $\partial\Omega$ is regular enough to allow an approximation of the domain $\Omega$ by triangular finite element mesh: Let  ${\mathscr T}_h = \{K_j\}_{j=1}^{N_t}$ be a conforming tessellation in triangles $K_j$, where $h$ is a symbolic mesh size parameter, and let $\Omega_h$ be a polygonal approximation of the domain,
\[
 \overline \Omega_h = \bigcup_{j=1}^{N_t}\overline K_j, 
\] 
Denoting by $\{\psi_j\}_{j=1}^{N_v}$ a set of appropriately chosen Lagrange basis functions related to the tessellation ${\mathscr T}_h$, where $N_v$ is the number of vertices of the mesh, we approximate the unknown $u$ by
\begin{equation}\label{u discr}
 u(x) \approx \widetilde u_h(x) = \sum_{j=1}^{N_v} u_j \psi_j(x), \quad u_j = u(x_j).
\end{equation}
Furthermore, let us write an approximate discrete model as
\[
 b = F(u_h, {\mathscr T}_h) + \varepsilon_h, \quad u_h = \left[\begin{array}{c} u_1 \\ \vdots \\ u_{N_v}\end{array}\right],
\] 
where  the error $\varepsilon_h$ contains the measurement error and the approximation error,
\begin{eqnarray*}
  b &=& F^*(u) + \varepsilon \\
    &=& F(u_h, {\mathscr T}_h)  + \underbrace{ \left[F^*(u) - F(u_h, {\mathscr T}_h) \right] + \varepsilon}_{=\varepsilon_h}.
\end{eqnarray*} 
While the modeling error cannot be estimated without knowning $u$, the Bayesian framework provides a natural solution. In the Bayesian approach, since $u$ is an unknown, it is modeled as a random variable $U$ with a prior probability distribution $\mu_U$. By denoting
\[
 {\mathscr F}: u \mapsto  m = \left[F^*(u) - F(u_h, {\mathscr T}_h) \right],
\]
we may define the probability distribution of the discrepancy, modeled as a random variable $M$, as the push-forward of the prior distribution,
\[
 \mu_M = {\mathscr F}_\# \mu_U.
\]
An analytic form of this distribution is usually not available, however, a computational approximation can be carried out using sampling techniques. In practice, instead of the accurate model $F^*$,  an ``accurate enough" model is used, and the modeling error sample is generated off-line, akin to network training in machine learning \cite{Arridge,CDSS,KSmoderror,Nissinen}.

To minimize the effect of the model discrepancy, one could use a mesh fine enough so that $\varepsilon_h$ is negligible, or, alternatively, find a way to estimate the modeling error distribution, which also requires an accurate model $F^*$, or a good approximation of it, to generate a model discrepancy sample. One of the main goals of this paper is to bypass the need for a globally accurate model by developing a computationally feasible way to perform selective mesh refinement. In other words, instead of using a costly dense mesh over the whole domain $\Omega$ to emulate $F^*$, the mesh is refined only where it is necessary to represent accurately the large gradients of the unknown. The same challenges are encountered in a posteriori mesh refining problems in standard finite element analysis, and the approach that we propose has been inspired largely by the contributions to finite element error analysis. In the adaptive meshing algorithms for approximating the solutions of second order PDEs of advection-diffusion-reaction (ADR) type by finite elements \cite{ZZ,MP,MPF}, referred to as the  a posteriori error estimate, an important role is played by the approximation error functional
\begin{equation}\label{post err}
 {\mathscr E}_h =\| \nabla u - \nabla u_h\|^2 \approx \| \nabla^*u - u_h\|^2,
\end{equation}
where $u$ is the true solution of the PDE, $u_h$ is the current Galerkin FEM approximation of the weak solution, and  $\nabla^* u$ is a computational proxy for the true solution, referred to as the  {\em recovered gradient}, a Cl\'{e}ment-type interpolant of a given order \cite{MP,MPF}. The adaptive meshing of the computational domain is guided by the  equi-distribution principle of the error: Each element in the new mesh should contribute approximately the same amount to the total error. The local error estimates can be expressed in terms of a metric that is related to the underlying meshing. The metric-based discretization is the main source of inspiration of the current approach \cite{Du,Hecht}.

Assume that $\mG = (g_{ij})$ is a metric tensor on $\Omega$, defined so that the length of a parametrized curve $x = \gamma(t)$ is given by the integral of the line element,
\[
 ds = \big((\gamma'(t))^\mT \mG \gamma'(t)\big)^{1/2} dt.
\]   
The principle guiding our meshing strategy is that the triangles $K_j$ have approximately the same diameter with respect to the metric distance for a given metric. Conversely, given a tessellation of $\Omega$ with triangles of different Euclidian size, we may define approximately a metric so that the elements have the same diameter in this metric. The connection between the tessellation and metric is developed in detail in Section~\ref{sec:metric}. Observe that the metric may be anisotropic, leading to triangles of varying eccentricity, defined in terms of the eccentricity of the Steiner ellipses circumscribing the triangle.

To design a metric that adapts to the unknown distributed parameter $u$, we would like to require that:
\begin{enumerate}
\item When $|\nabla u|$ is large, the triangles defined by the metric are relatively small in the Euclidian metric;
\item For anisotropic triangles, the shorter semi-axes of the Steiner ellipses are parallel to the gradient;
\item For $|\nabla u|$ small, the metric is isotropic, corresponding to spherical Steiner ellipses.
\end{enumerate}

A candidate metric satisfying these requirements can be written as
\[
 \mG(x) = \theta_{\|}(x) e_{\|}(x)e_{\|}(x)^\mT + \theta_{\perp}(x) e_\perp(x) e_\perp(x)^\mT,
\]
where $(e_{\|}(x),e_\perp(x))$ is a local orthonormal frame such that $e_{\|}(x)$ is parallel to the gradient of $u(x)$, and
\[
 \theta_{\|}(x), \; \theta_\perp(x)\geq  \theta_{\rm min}>0
\]
are non-decreasing scalar functions of the norm of the gradient of $u(x)$, with the property that for $|\nabla u(x)|$ large,
\[
 \frac{\theta_{\|}(x)}{\theta_\perp(x)} = \alpha \geq 1,
\] 
where $\alpha$ is the {\em maximum anisotropy ratio}, and for $|\nabla u(x)|$ small, ${\theta_{\|}(x)} = {\theta_\perp(x)}$.
 Thus, for $|\nabla u(x)|$ large, the triangles are small, with width in the direction of $\nabla u(x)$ by a factor $\sqrt{\alpha}$ smaller than the the width in the orthogonal direction, while for small gradients, the triangles are relatively large and isotropic. The exact  definitions of the functions $\theta_{\|}$ and $\theta_\perp$ and of  the metric respecting these principles are given in Section~\ref{sec:metric}.

Given a metric and the corresponding tessellation $({\mathscr T}_h)$, the discrete inverse problem to be solved corresponds to the model
\begin{equation}\label{model Mh}
 {\mathscr M}_h : \quad b = F\big(u_h,{\mathscr T}_h\big) + \varepsilon_h.
\end{equation}
When solving the inverse problem in the Bayesian framework, it is necessary to specify in which way the quality of the solution will be assessed, to define a feedback mechanism to update the metric, and to provide a way to implement  the mesh refinement. This feedback is defined by using hierarchical Bayesian models, and therefore implemented as part of the prior for $u$. In this manner, {\em finding the optimal metric, and therefore the optimal discretization, becomes part of the inverse problem itself.}

The reconstruction problem will be addressed in a hierarchical Bayesian framework where the prior assumption that $u$ is piecewise constant, or nearly piecewise constant, is tantamount to assuming that the discretized gradient can be represented as a sparse or compressible vector. To this end, we set up a Bayesian model that favors sparse solutions of the discretized gradient of $u$.
We want to develop an algorithm suitable for estimating a distributed variable $u$ with the a priori belief that the numerical approximation of the total variation of $u$ is limited. We begin by setting up a prior expressing the belief that the discrete approximation of the gradient  of $\widetilde u_h$ given by (\ref{u discr}) is compressible. We begin with representing $\nabla \widetilde u_h$ in terms of a Whitney basis,
\begin{equation}\label{nabla u}
 \nabla \widetilde u_h(x) = \sum_{\ell=1}^{N_e} z_\ell w_\ell(x),
\end{equation}
where $N_e$ is the number of edges in the triangular mesh generated by the tessellation ${\mathscr T}_h$, and $w_\ell$ is the Whitney basis function associated to the $\ell$th edge. Denoting the first order piecewise linear nodal basis functions over the current mesh by $\psi_j$, the Whitney basis function associated with a directed  edge $e_\ell = \{v_j,v_k\}$ is given by
\[
 w_\ell = \psi_j\nabla \psi_k - \psi_k\nabla\psi_j,
\]
the gradients being understood in the weak sense. The relation between the coefficients $z_\ell$ and the nodal values of the function $u$ is discussed in detail in the next section. With the representation introduced above, we want to find $u$ so that the coefficient vector $z$ is sparse, or more generally, compressible.
  
To set up the hierarchical model first in the abstract setting, consider a discrete inverse problem,
\begin{equation}\label{forward}
 B = f(Z) + E,
\end{equation}
where $Z$ and $E$ are random variables taking values in $\R^n$ and $\R^m$ representing the unknown and noise, respectively, and $B$ is the $\R^m$-valued random variable representing the data. Moreover, we assume that the random variable $Z$ is sparse or compressible, meaning that only a limited number of its components are above some small threshold value.
We express our information about $Z$ in terms of a prior density conditional on a hyperparameter $\Theta$, denoted by $\pi_{Z\mid\Theta}(z\mid \theta)$. In turn, the hyperparameter itself is a random variable in $\R^k$, and we define a hyperprior density $\pi_\Theta(\theta)$ to describe the prior belief about it. The joint prior density of the pair $(Z,\Theta)$ is given by
\[
 \pi_{Z,\Theta}(z,\theta) = \pi_{Z\mid\Theta}(z\mid\theta)\pi_\Theta(\theta).
\]  
Next we write the likelihood density $\pi_{B\mid Z}(b\mid z)$ based on the forward model (\ref{forward}).  If $Z$ and $E$ are assumed mutually independent, the likelihood is of the form
\begin{equation}\label{likelihood}
 \pi_{B\mid Z}(b\mid z) = \pi_E\big( b - f(z)\big),
\end{equation}
where $\pi_E$ is the probability density of the additive noise. The posterior density is given by Bayes' formula,
\begin{equation}\label{posterior}
 \pi_{Z,\Theta}(z,\theta\mid b) \propto  \pi_{Z,\Theta}(z,\theta) \pi_{B\mid Z}(b\mid z) = 
 \pi_{Z\mid\Theta}(z\mid\theta)\pi_\Theta(\theta) \pi_{B\mid Z}(b\mid z).
\end{equation}
Given the posterior density, it is possible to compute the Maximum A Posteriori (MAP) estimate $(z_{\rm MAP},\theta_{\rm MAP})$, a maximizer of the posterior density above, assuming that such maximizer exists.
In a series of articles \cite{MEG,CSS,CPPSV,CPSS}, the authors have developed a methodology based on a family of hierarchical models that promote sparsity of the variable $z_{\rm MAP}$ by favoring solutions for which the majority of the components of $z$ are either zero or of negligible size. Since the algorithm proposed in this work for the adaptive remeshing is based on those ideas, for completeness we provide a more in-depth review of the sparsity promoting hierarchical prior methodology in Section~\ref{sec:IAS}, while only presenting a general outline here.  

Let $Z$ be a compressible random variable, with zero mean hierarchical conditionally Gaussian prior density
\[
 \pi_{Z\mid\Theta}(z\mid\theta) = \prod_{j=1}^n \frac{1}{\sqrt{2\pi\theta_j}}{\rm exp}\left(-\frac{z_j^2}{2\theta_j}\right) \propto {\rm exp}\left(-\frac 12 \sum_{j=1}^n\frac{z_j^2}{\theta_j} - \frac 12 \sum_{j=1}^n\log \theta_j\right).
\]  
When the prior variance $\theta_j$ is small, it is expected that the variable $Z_j$ takes on a small value, too. Since part of the the problem is to identify the few components of $Z$ that are large, the hypermodel for the variance $\Theta$ should favor small values most of the time, while allowing occasional large outliers. Among the several fat-tailed distributions that could be chosen for the he components of $\Theta$, for reasons of computational convenience we choose the generalized gamma distribution, i.e., 
\[
 \Theta_j \sim \pi_{\Theta_j}(\theta_j \mid \vartheta_j,\beta,r)  = \frac{|r|}{\Gamma(\beta)\vartheta_j}\left(\frac{\theta_j}{\vartheta_j}\right)^{r\beta -1}
 {\rm exp}\left( - \left(\frac{\theta_j}{\vartheta_j}\right)^r\right), \quad 1\leq j\leq n,
\]
where $r\neq 0$, and the shape parameter $\beta>0$ is assumed to be the same for all $j$ while the scale parameter $\vartheta_j$ may differ for every $j$.  Assuming that the variables $\Theta_j$ are mutually independent,  the joint prior model then becomes
\begin{eqnarray}\label{joint prior}
 \pi_{Z,\Theta}(z,\theta) &=& \pi_{Z\mid\Theta}(z\mid\theta)\pi_{\Theta}(\theta\mid \vartheta,\beta,r) \nonumber \\
 &\propto& {\rm exp}\left( -\frac 12 \sum_{j=1}^n \frac{z_j^2}{\theta_j}  - \sum_{j=1}^n\left(\frac{\theta_j}{\vartheta_j}\right)^r +\left(r\beta - \frac 32\right) \sum_{j=1}^n \log\frac{\theta_j}{\vartheta_j} \right).
 \end{eqnarray}

In the case  where the additive noise is Gaussian, $E\sim{\mathcal N}(0,\mSigma)$, where $\mSigma\in\R^{m\times m}$ is a symmetric positive definite covariance matrix, the likelihood density (\ref{likelihood}) is the form
\[
 \pi_{B\mid Z}(b\mid z) \propto {\rm exp}\left(-\frac 12 \big(b - f(z)\big)^\mT\mSigma^{-1}\big(b - f(z)\big)\right),
\]
and the expression for the posterior density  (\ref{posterior}) becomes
\begin{eqnarray*}
 && \pi_{Z,\Theta\mid B}(z,\theta\mid b) \\
 && \phantom{XX} \propto {\rm exp}\left(-\frac 12 \big(b - f(z)\big)^\mT\mSigma^{-1}\big(b - f(z)\big)
  -\frac 12 \sum_{j=1}^n \frac{z_j^2}{\theta_j}  - \sum_{j=1}^n\left(\frac{\theta_j}{\vartheta_j}\right)^r +\left(r\beta - \frac 32\right) \sum_{j=1}^n \log\frac{\theta_j}{\vartheta_j} 
  \right).
\end{eqnarray*}
Algorithms for finding a sparse or compressible estimate maximizing the above expression, or, equivalently, minimizing the Gibbs energy,
\begin{equation}\label{Gibbs}
{\mathscr E}(z,\theta) =
\frac 12 \big(b - f(z)\big)^\mT\mSigma^{-1}\big(b - f(z)\big)+
  \frac 12 \sum_{j=1}^n \frac{z_j^2}{\theta_j}  + \sum_{j=1}^n\left(\frac{\theta_j}{\vartheta_j}\right)^r -\left(r\beta - \frac 32\right) \sum_{j=1}^n \log\frac{\theta_j}{\vartheta_j} ,
\end{equation}
are discussed in \cite{MEG,CPS,CPSS,CPPSV,CSS,PCS}. An overview of the methodology will be given in Section~\ref{sec:IAS}.

The Bayesian adaptive (anisotropic) remeshing algorithms that we propose here can be summarized as follows.
\begin{enumerate}
\item Given the current tessellation ${\mathscr T}_h$,
\begin{enumerate}
\item Discretize the forward problem (\ref{model Mh}), representing the unknown $\widetilde u_h$ in terms of the degrees of freedom $z_j$ corresponding to the gradient;
\item Compute the MAP estimate $(z_{\rm MAP},\theta_{\rm MAP})$ using the hierarchical Bayesian model;
\end{enumerate}
\item Interpolate the MAP estimate for $\nabla \widetilde u_h$ and update the metric;
\item Update the tessellation ${\mathscr T}_h$ in accordance with the updated metric;
\item If convergence criterion is met, stop, otherwise repeat from 1. 
\end{enumerate}
Several important details are omitted in the outline above, including how the the gradient is interpolated to allow point-wise estimates, and metric is related to the interpolated gradient.  These and other details are presented below in separate subsections.  The performance of the algorithm is then illustrated with computed examples. 

\section{Detailed derivation of the Bayesian adaptive meshing algorithm}

In this section, we discuss the details of the various steps outlined above, and we conclude it with a detailed algorithm. We begin with the details of the model discretization.

\subsection{Model discretization}

Let ${\mathscr T}_h$ denote the current tessellation of the domain $\Omega_h\approx \Omega$, which for simplicity we assume to be simply connected and denote by $\{v_j\}_{j=1}^{N_v}$ the vertices of the mesh. We assume that (\ref{u discr}) is a piecewise linear approximation of the unknown function $u$, that is, the basis functions $\{\psi_j\}_{j=1}^{N_v}$ are the standard piecewise linear  Lagrange basis functions,
\[
 \psi_j(v_k) = \delta_{jk}.
\]
We assume that the  first $n_v$ nodes are interior nodes and that the values of $u$ at the boundary nodes are known; without loss of generality, we assume that
\[
 u(v_j) = 0, \quad n_v+1\leq j\leq N_v.
\] 
The number of boundary nodes is denoted by $n_v^b = N_v - n_v$.

\begin{figure}
\centerline{
\includegraphics[width=10cm]{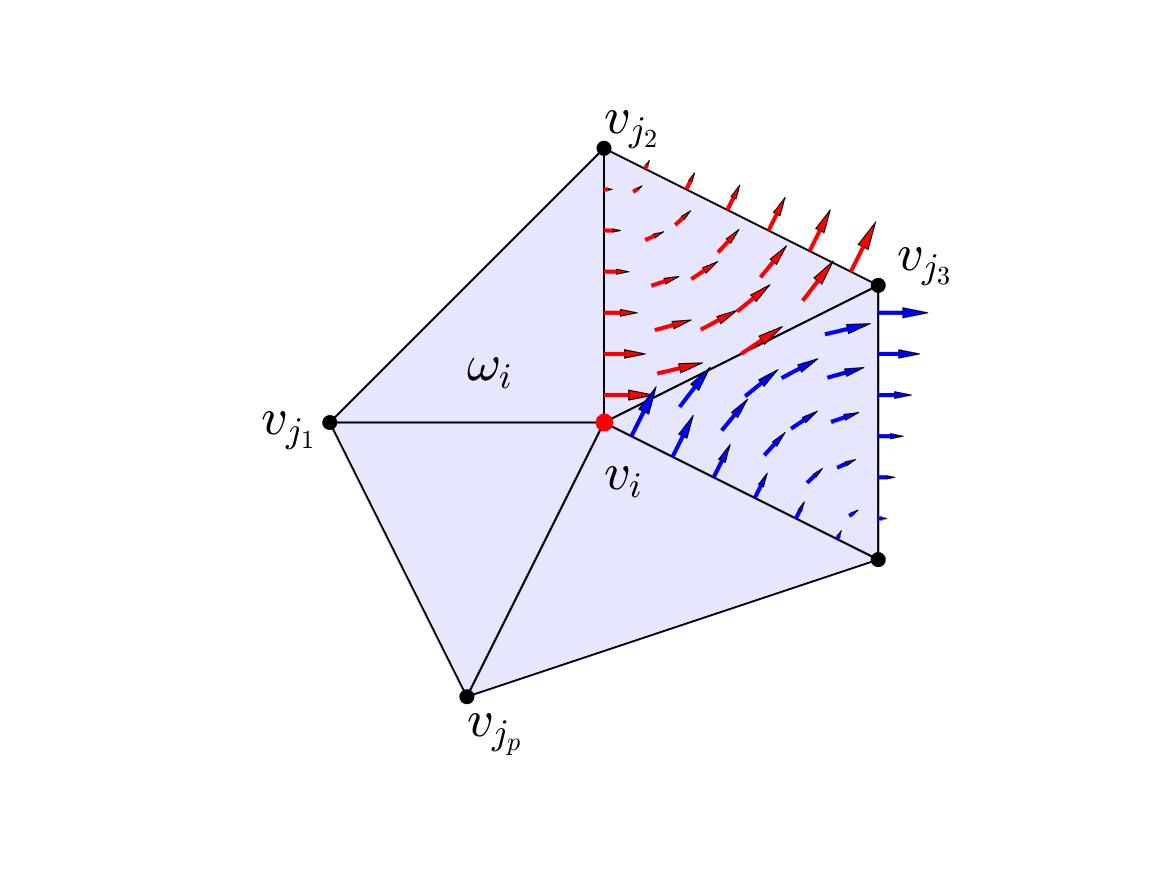}
}
\caption{\label{fig:patch Whitney} The local patch $\omega_i$ associated to the vertex $v_i$. In the figure, the Whitney basis function associated to the edge $\{v_i,v_{j_3}\}$ is shown as a vector field plot. Outside the two triangles sharing the edge the basis function vanishes. Observe that the basis function is discontinuous, while having continuous tangential component across the edges.}
\end{figure}
Consider the edge elements of Whitney type: If $e_\ell = \{v_j,v_k\}$ is the directed edge connecting the nodes $v_j$ and $v_k$, we define 
\[
 w_\ell = \psi_j\nabla\psi_k - \psi_k\nabla\psi_j, \quad 1\leq \ell\leq  N_e,
\]
the gradients understood as weak derivatives of the piecewise linear basis functions; see Figure~\ref{fig:patch Whitney}.  Given the approximation (\ref{u discr}) in terms of the nodal values $u_j$, the passage to the approximation (\ref{nabla u}) is straightforward, and can be found in the literature \cite{Bossavit,Bossavit2}, however, for completeness, we derive it in detail below. Let $v_j$ be an arbitrary node, and denote its neighboring nodes by $v_{j_1}, \ldots,v_{j_p}$.  Let $I^i$ denote the set of indices of the edges with the vertex $v_i$ as one of the endpoints.
In the patch neighborhood $\omega_i$, we have
\[
 \sum_{\ell = 1}^p \psi_{j_\ell}(x) + \psi_i(x) = 1,\quad x\in\omega_i,
\]
implying that
\[
 \nabla \psi_i(x) = -\sum_{\ell = 1}^p \nabla \psi_{j_\ell}(x).
\]  
Therefore,
\begin{eqnarray*}
\sum_{\ell = 1}^p\big(\psi_i(x) \nabla\psi_{j_\ell}(x) - \psi_{j_\ell}(x) \nabla\psi_i(x) \big)  
&=& \psi_i(x)\sum_{\ell=1}^p \nabla\psi_{j_\ell}(x)  - \nabla\psi_i(x)\sum_{j=1}^p \psi_{j_\ell}(x) \\
&=& -\psi_i(x) \nabla\psi_i(x)  -\nabla \psi_i(x)\big(1 - \psi_i(x)\big)  \\
&=& -\nabla\psi_i(x),
\end{eqnarray*}
that is,
\[
 \nabla \psi_i(x) = -\sum_{\mu \in I^i}  s_\mu^i w_\mu(x),
\] 
where $s_\mu^i = \pm 1$, the sign depending on whether $v_i$ is the tail or the head of the edge $e_\mu$.  This shows that the gradients of the basis functions are in the span of the Whitney elements. Moreover, it follows from this representation that
\[
 \nabla \widetilde u_h(x) =  \sum_{i = 1}^{N_v} u_i\nabla\psi_i(x) = -  \sum_{i = 1}^{N_v}  u_i \sum_{\mu \in I^i}  s_\mu^i w_\mu(x).
\]  
Every edge $\mu$ appears exactly twice in the sum, corresponding to the end edges, and if $e_\mu = \{v_j,v_k\}$, we have $s^j_\mu = -1$ and $s^k_\mu = 1$. Therefore 
\[
  \nabla \widetilde u_h(x)  = \sum_{\ell=1}^{N_e} z_\ell w_\ell(x), \quad z_\ell = u_j - u_k.
\]  
Hence, the coefficients $z_\ell$ are obtained by a finite difference formula that is independent of the geometry but depends only on the topology of the mesh. We write the relation in matrix form as
\[
 z = \mL_{\rm full}\, u, \quad \mL_{\rm full}\in\R^{N_e\times N_v},
\]
where each row of the matrix $\mL_{\rm full}$ contains exactly two non-zero entries, $\big(\mL_{\rm full}\big)_{\ell j} = 1$, $\big(\mL_{\rm full}\big)_{\ell k} = -1$. 

To modify the formula so that it contains only the nodal values in the interior nodes, we observe first that the matrix $\mL_{\rm full}$ has a one-dimensional null space,
\begin{equation}\label{null}
 {\mathcal N}\big(\mL_{\rm full}\big) = {\rm span}\{ {\bf 1}\}, \quad {\bf 1} = \left[\begin{array}{c} 1 \\ \vdots \\ 1\end{array}\right]\in\R^{N_v}.
\end{equation}
Next we partition the matrix into the columns corresponding to interior and those corresponding to boundary nodes,
\[
 \mL_{\rm full} = \left[\begin{array}{cc} \mL_{\rm int} & \mL_{\rm bdry}\end{array}\right], \quad \mL_{\rm int} \in\R^{N_e\times n_v},
\] 
To show that $\mL_{\rm int}$ has a trivial null space, let $z_0\in {\mathcal N}(\mL_{\rm int})$. Then,
\[
 \mL_{\rm full}\left[\begin{array}{c} z_0 \\ 0\end{array}\right] = \mL_{\rm int} \,z_0 + \mL_{\rm bdry}\,0 = 0,
\]
and from (\ref{null}) it follows that $z_0=0$. Finally, after permuting the rows so that those corresponding to the boundary edges are at the bottom, we partition $\mL_{\rm int}$ as
\[
 \mL_{\rm int} = \left[\begin{array}{cc} \mL \\ \mL'\end{array}\right], \quad \mL \in\R^{n_e\times n_v}.
\]
Since for every $x\in\R^{n_v}$, we must have
\[
 \mL_{\rm int}\, x = \left[\begin{array}{c} \mL x \\ \mL' x\end{array}\right] =  \left[\begin{array}{c} \mL x \\ 0\end{array}\right] ,
\]
we conclude that $\mL' = \mO_{(N_e-n_e)\times n_v}$.  In the following, we denote by $z\in\R^{n_e}$ the increments over edges that have at least one end point in the interior, thus excluding the zero increments over boundary edges, and write
\begin{equation}\label{z=Lu}
 z = \mL u, \quad \mL\in\R^{n_e\times n_v},
\end{equation}
where the matrix is full rank, as $n_e = n_v + N_t -1\geq n_v$ by the definition of the Euler characteristics. 

\subsection{Tessellation and metric}\label{sec:metric}

\begin{figure}
\centerline{
\includegraphics[width=9cm]{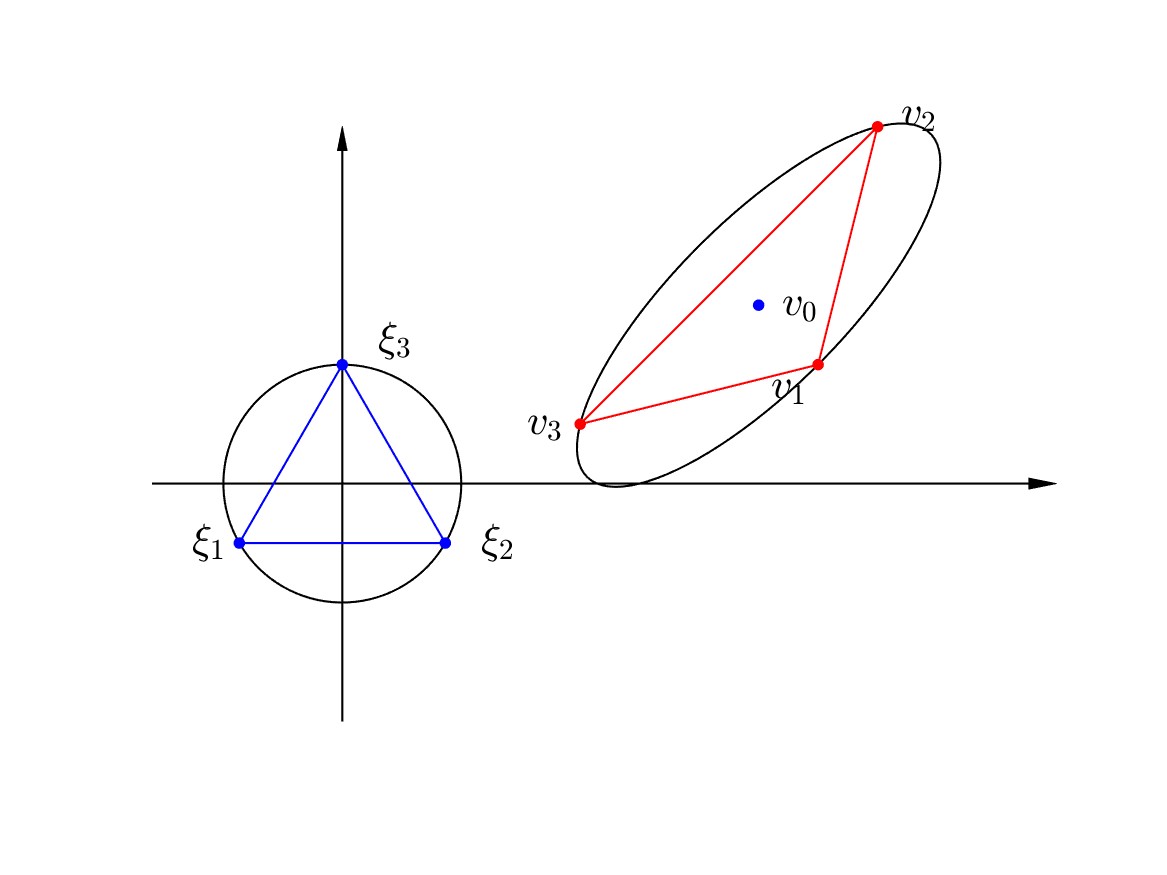}
\includegraphics[width=9cm]{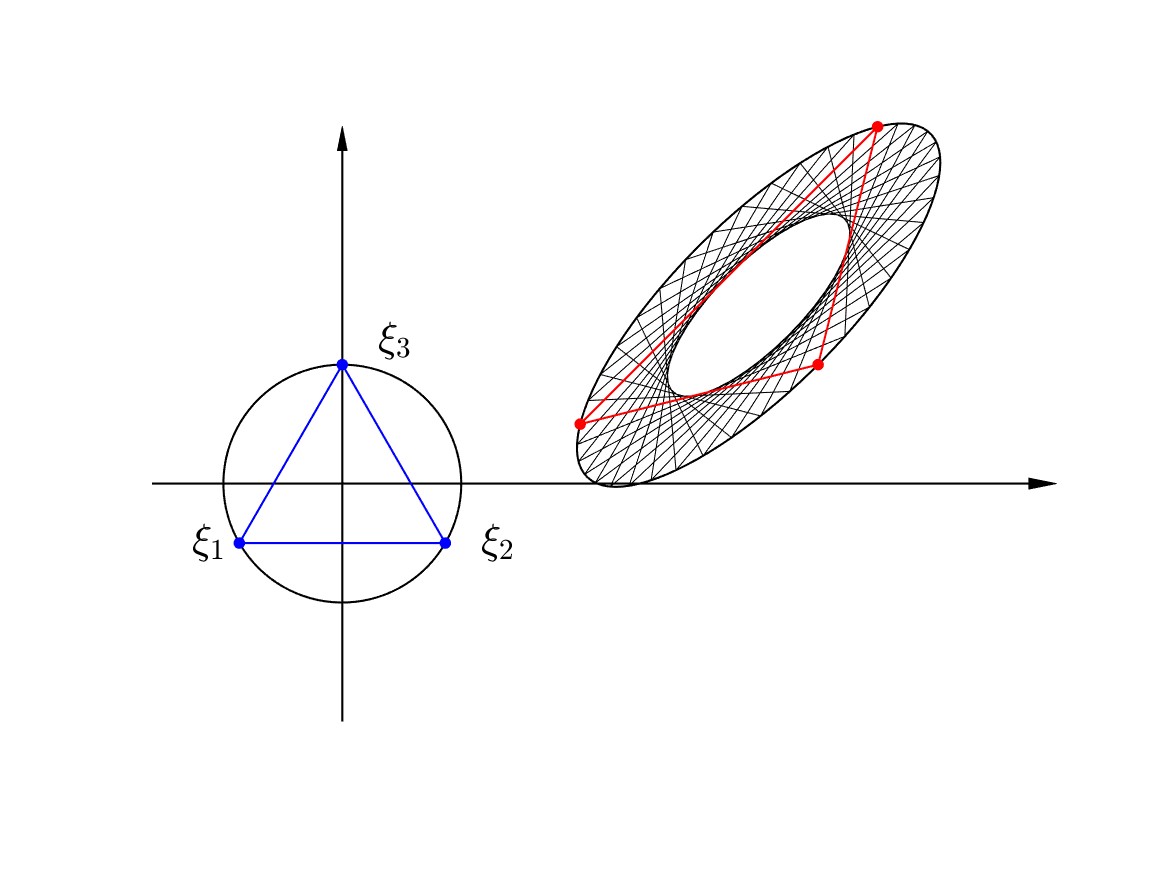}
}
\caption{\label{fig:ellipses} Left panel: An equilateral reference triangle $\widehat K$  plotted in blue, and a physical triangle $K$ plotted in red, along with the Steiner circumellipse. On the right, a family of triangles with the same polar matrix $\mP$ but with different rotations $\mW$, all sharing the same Steiner ellipse and thus having the same eccentricity.}
\end{figure}

We start with a brief summary of the metric-based tessellation algorithm.  Consider a non-isoparametric triangle $K\in{\mathscr T}_h$ with vertices $v_1$, $v_2$ and $v_3$ ordered counterclockwise, and let $\widehat K$ denote an equilateral reference triangle mapped onto $K$ by the affine mapping
\[ 
 F_K : \widehat K \to K,\quad \xi \mapsto v_0 + \mF \xi ,\quad  v_0 = \frac 13\sum_{j=1}^3 v_j. 
 \]
Consider the polar decomposition of $\mF$, that can be derived from the singular value decomposition by letting 
\[
 \mF = \mU \mD \mV^\mT = \big(\mU \mD \mU^\mT\big) \big(\mU \mV^\mT) = \mP \mW, 
\]
where $\mW$ is an orthogonal matrix and $\mP$ is a symmetric positive definite (SPD) matrix. The formula defining the Steiner circumellipse is 
\[
  s\mapsto v_0 + x(s)^\mT\mP x(s), \quad x(s) = \left[\begin{array}{c} \cos 2\pi s \\ \sin 2\pi s\end{array}\right], \quad 0\leq s\leq 1,
\]   
and, conversely, the family of triangles with the same Steiner ellipse is given by
\[
 \big\{v_0 + \mP\mW(\widehat K)\mid \mW \in {\rm SO}(2)\big\},
\]
see Figure~\ref{fig:ellipses}. 

Given the tessellation, we define the piecewise constant metric over the approximation $\Omega_h$ of the domain $\Omega$ by setting
\[
 \mG\big|_{K} = \mP^{-2}, \quad K\in{\mathscr T}_h.
\]
The arc length element of a line segment across the triangle $K$, given by $\gamma(t) = v_0 + \mF \xi(t)$, where  $\xi(t)\in\widehat K$ and $\mF = \mP\mW$ is then given by
\[
 ds =  \big(\gamma'(t)^\mT\mG  \gamma'(t)\big)^{1/2} dt  
 =  \big(\xi'(t)^\mT\mW ^\mT \mP\mG  \mP\mW \xi'(t)\big)^{1/2} dt 
= |\xi'(t)| dt,
\]
implying that with respect to this metric, the triangles are of equal size. Conversely, we consider the following meshing principle: Given a metric $\mG=\mG(x)$, $x\in\Omega$, find a tessellation such that the triangles have approximately the same size with respect to the metric. The problem can be set up as follows.

\begin{problem}
Given a metric $\mG$ in $\Omega$ and $0<h_{\rm min}<h_{\rm max}$, find a tessellation $\{K_j\}$ of $\Omega$ such that 
$h_{\rm min}\leq {\rm diam}(K_j)\leq h_{\rm max}$, minimizing the discrepancy
\[
 \max_j \|\mG\big|_{K_j} - \mP_{K_j}^{-2}\|_\infty
\] 
over the tessellation.
\end{problem}

Observe that the bounds $h_{\rm min}$ and $h_{\rm max}$ are safeguards that could be built in the metric itself. In this work, the mesh generation is done by employing an existing software FreeFEM and particularly the function {\tt adaptmesh}, see \cite{Hecht} for details.

\subsection{Solving the MAP estimate: IAS algorithm}\label{sec:IAS}

Consider the Gibbs energy (\ref{Gibbs}) in the case when the forward model corresponds to a linear observation model, $f(z) = \mA z$. Let $u\in\R^{n_v}$ denote the nodal values of the discretized function of interest, and let $z\in \R^{n_e}$ denote the coefficients in (\ref{nabla u}) of the gradient in terms of the Whitney basis, related to $u$ through the equation (\ref{z=Lu}).  Since our goal is is to write the sparsity-promoting prior for $z$, it is necessary to express the likelihood in terms of $z$.

Consider the QR-decomposition of the matrix $\mL$,
\[
 \mL = \mQ \mR = \left[\begin{array}{cc} \mQ_1 & \mQ_2\end{array}\right]\left[\begin{array}{c} \mR_1 \\ \mO\end{array}\right],
\]
where $\mR_1\in\R^{n_v\times n_v}$ is an upper triangular matrix, $\mQ\in\R^{n_e\times n_e}$ is orthogonal, $\mO\in\R^{(n_e-n_v)\times n_v}$ is a null matrix, and $\mQ_1$ and $\mQ_2$ are a partitioning of $\mQ$.  
Since $\mL$ has rank $n_v$, the matrix $\mR_1$ is invertible,  hence (\ref{z=Lu}) is equivalent to 
\[
 u = \mR_1^{-1}\mQ_1^\mT z, \quad \mQ_2^\mT z = 0,
\]
where the second equation corresponds to the compatibility condition of vanishing circulations over each triengle, expressed in the orthonormal basis given by $\mQ$. Introducing
\[
\mA_1 = \mA \mR_1^{-1} \mQ_1^\mT,
\]
we may write the Gibbs energy in terms of $z$ as
\begin{equation}\label{gibbs z}
{\mathscr E}(z,\theta) =  \frac 12 \big\|\mS\big(b - \mA_1 z\big)\big\|^2 +
  \frac 12 \sum_{j=1}^{n_e} \frac{z_j^2}{\theta_j}  + \sum_{j=1}^{n_e}\left(\frac{\theta_j}{\vartheta_j}\right)^r -\left(r\beta - \frac 32\right) \sum_{j=1}^{n_e} \log\frac{\theta_j}{\vartheta_j},\quad \mbox{with $\mQ_2^\mT z = 0$.}
\end{equation}
where $\mSigma^{-1} = \mS^\mT \mS$ is a symmetric factorization of the precision matrix of the noise. Without loss of generality, we may assume that $\mS = \mI$, that is, the additive noise is whitened by scaling $\mA_1$ and $b$ by the factor $\mS$. 

The basic IAS algorithm seeks a minimizer of the Gibbs energy by alternating the following minimization steps:
\begin{enumerate}
\item Given the current approximation of $\theta$, update $z$ by minimizing
\[
{\mathscr E}_\theta(z) =  \frac 12 \big\|\mS\big(b - \mA_1 z\big)\big\|^2 +
  \frac 12 \sum_{j=1}^{n_e} \frac{z_j^2}{\theta_j} , \quad \mbox{with  $\mQ_2^\mT z = 0$.} 
\]
\item Given the current approximation of $z$, update $\theta$ by minimizing
\[
  {\mathscr E}_z(\theta) = 
  \frac 12 \sum_{j=1}^{n_e} \frac{z_j^2}{\theta_j}  + \sum_{j=1}^{n_e}\left(\frac{\theta_j}{\vartheta_j}\right)^r -\left(r\beta - \frac 32\right) \sum_{j=1}^{n_e} \log\frac{\theta_j}{\vartheta_j}.
 \] 
 \end{enumerate}
We remark that the updating of $z$ is tantamount to solving a least squares problem, and the updating of $\theta$ can be done componentwise. Indeed, by differentiating with respect to $\theta_j$, the first order optimality condition yields
\begin{equation}\label{theta cond}
  \frac{\partial{\mathscr E}_z}{\partial \theta_j}(\theta) = 
  -\frac 12 \frac{z_j^2}{\theta_j^2}  + r \frac 1{\theta_j}\left(\frac{\theta_j}{\vartheta_j}\right)^r -\frac{1}{\theta_j}\left(r\beta - \frac 32\right) =0.
 \end{equation}
In the special cases $r = \pm 1$ this equation admits an explicit solution, while in general, a numerical scheme is required, see \cite{CPSS}.  

While in principle rather straightforward, there are several details in the algorithm that need to be addressed properly. Before discussing the fine tuning of the algorithm, we will look at the hyperparameters.

\subsubsection{Hyperparameters and sensitivity scaling}

Let's start by considering the hypermodel $r=1$ corresponding to the gamma distribution. In the literature \cite{MEG,CSS,CSbook}, the following results have been proved.

\begin{theorem}\label{th:gamma}
\begin{itemize}
\item[(a)] The Gibbs energy (\ref{gibbs z}) for the gamma hyperprior ($r=1$) is a convex functional of $(x,\theta)$ with a unique minimizer $(\widehat z,\widehat\theta)$, and the IAS iteration converges to this minimizer.
\item[(b)] Moreover, when $\eta = \beta -3/2>0$ converges to zero, the minimizer $\widehat z$ converges to the minimizer of the functional
\begin{equation}\label{ell 1}
 {\mathscr E}_0(z) = \frac 12 \big\|b - \mA_1 z\|^2 + \sqrt{2} \sum_{j=1}^{n_e} \frac{ |z_j|}{\sqrt{\vartheta_j}}.
\end{equation} 
\end{itemize}
\end{theorem}
The second part of the above theorem suggests that to guarantee sparsity of the solution the variable $\eta$ should be set to a  small value. 
Observe that in (\ref{ell 1}) there is no explicit dependence on the noise level, because of the whitening assumption. The formula indicates that the choice of the values of the hyperparameters $\vartheta_j$ must be related to the sensitivity of the data to the components $z_j$. Indeed, it is well understood that in inverse problems, the weighting of the components in the penalty term is essential when the sensitivity of the data to different components varies: Without sensitivity weighting, optimization algorithms will converge to solutions that explain the data in terms of the variables to which the data are most sensitive,  see  \cite{CSbook} for a  discussion and, e.g.,  \cite{Li1,Li2,Lin}  in the context of geophysics and medical imaging.  In the optimization literature, sensitivity weighting is referred to as scaling of the problem \cite{DennisSchnabel}. This observation may appear to be in conflict with the Bayesian philosophy according to which the prior is independent of the observation model. {In \cite{CPPSV,CPSS},  the authors established a relationship between the hyperparameters and the prior information about  signal-to-noise ratio (SNR), thus providing a consistent interpretation of sensitivity weighting. In the same article, the concept of {\em SNR-exchangeability} was introduced: A prior model satisfies the SNR-exchangeability condition, if any subset of components having the same cardinality can lead to the expected SNR. In other words, the SNR-exchangeability principle states that all subsets of a given cardinality of the coordinates of the unknown  must have the same possibility to explain the data within the presumably known SNR level. It was shown also that if there is SNR-exchangeability, the hyperparameters $\vartheta_j$ must satisfy
\begin{equation}\label{sens scaling}
 \vartheta_j = \vartheta^*_r \| \mA_1 e_j\|^2,
 \end{equation}
 i.e., $\vartheta_j$ is proportional to the squared norm of the $j$th column of the matrix $\mA_1$. Interestingly, in light of formula (\ref{ell 1}), this principle leads exactly to the classical sensitivity weighting or scaling referred to above. The parameter $\vartheta^*_r>0 $ depends on the hypermodel and the estimated SNR; we refer to the original articles \cite{CPPSV,CPSS} for the details.

\subsubsection{Non-convexity and hybrid schemes}

The discussion in the previous subsection was restricted to the IAS algorithm for $r=1$. In \cite{CPS,CPSS}, the IAS algorithm with $r<1$ was discussed in detail, and it was demonstrated that the Gibbs energy functional to be minimized is non-convex. In those cases the IAS algorithm converges rapidly, however, the minimizer found is typically a local minimizer that depends on the starting point. In order to take advantage of the fast convergence rate and greedier sparsity promotion of the IAS algorithm with $r<1$ while avoiding to stop at unsuitable local minimizers, a {\em hybrid IAS algorithm} was proposed. The hybrid approach consists of two phases: Phase I, using  the hypermodel $r=1$, lets the algorithm find a reasonable approximation of the unique minimizer of the Gibbs energy functional. Phase II switches to a hyperprior with $r<1$, and the IAS algorithm is run to convergence. The heuristic idea is that this way, the greedier Phase II converges to a minimizer which is close to the global minimizer of Phase I. 

To match the hyperparameters in Phases I and II so as to keep consistency throughout the entire process, the following matching conditions were proposed in \cite{sampler} (see also \cite{CPSS}):
\begin{enumerate}
\item When $z_j=0$, the formula (\ref{theta cond}) should yield the same value for $\theta_j$ with both hyperpriors,
\item The expected value of $\theta_j$ in Phase I and in Phase II should remain the same.
\end{enumerate}
Assuming that the parameters of Phase I are set, that is, $r=1$, and $\beta_1 = 3/2 + \eta$,  these two conditions lead to equations 
\begin{equation}\label{compatibility}
 \vartheta_1^*\eta =  \vartheta_2^*\left(\beta_2 - \frac{3}{2r_2}\right)^{1/r_2},
\end{equation}
and
\begin{equation}\label{compatibility2}
 \vartheta_1^*\left(\frac 32 + \eta\right)
= \vartheta_2^*\frac{\Gamma(\beta_2 + \frac 1{r_2})}{\Gamma(\beta_2)},
\end{equation} 
from which the pair $(\beta_2,\vartheta_2^*)$ can determined. Closed form solutions for the special cases $r = \pm 1/2$ and $r_2 = -1$ were given in \cite{sampler}. In our computed examples, we will use the hybrid model with $r_2 = 1/2$,  for which  
\[
 \beta_2  = \frac{6\mu+1  + \sqrt{48 \mu +1}}{2(\mu-1)}, \quad  \mbox{where } \mu = 1+\frac3{2\eta},
\]
and
\[
 \vartheta_2^* =  \vartheta_1^* \frac{\eta}{(\beta_2 -3)^2}.
\] 
With these choices of the hyperparameters, only two scalar variables,  $\eta$ and $\vartheta_1^*$, remain to be set. The discussion of the  choice of their values is postponed to the section of computed examples.

\subsubsection{Implementational details}

To guarantee that the IAS algorithm produces iterates of $z$ satisfying the compatibility condition of vanishing circulation around each element, we organize the computations as indicated in \cite{CCPS}. In the update of $z$, with the current value of $\theta$, we introduce  auxiliary variables
\[
 w_j = \frac {z_j}{\sqrt{\theta_j}},
\]
and write the objective function ${\mathscr E}_\theta(z)$ in terms of $w$.  Defining $\mD_\theta = {\rm diag}(\theta)$, we have
\[
 w = \mD_\theta^{-1/2} \mL u = \mL_\theta u, \quad \mL_\theta =  \mD_\theta^{-1/2} \mL.
\]
Next we compute the reduced QR factorization of $\mL_\theta$,
\[
 \mL_\theta = \mQ_\theta \mR_\theta, \quad \mQ_\theta \in \R^{n_e\times n_v}, \quad \mR_\theta \in \R^{n_v\times n_v},
\]
where $\mR_\theta$ is upper triangular and invertible, and write  

\[
  \frac 12 \big\|b - \mA u \big\|^2 +
  \frac 12 \sum_{j=1}^n \frac{z_j^2}{\theta_j}  =  \frac 12 \|b - \mA \mR_\theta^{-1}\mQ_\theta^\mT w\|^2 +\frac 12  \|w\|^2.
\]   
It was proved  in \cite{CCPS} that the minimizer of this expression satisfies automatically the compatibility condition. To compute efficiently an approximate minimizer, we observe that the least squares solution of the above functional is the standard Tikhonov regularized solution of the linear inverse problem 
\[
 b = \mA_\theta w, \quad \mA_\theta =  \mA \mR_\theta^{-1}\mQ_\theta^\mT ,
\] 
with the regularization parameter equal to unity. The Tikhonov regularized solution is the least squares solution of an overdetermined linear system, and finding it may be computationally very costly. We can reduce the computational complexity, approximating the Tikhonov regularized solution  by applying the Conjugate Gradient for Least Squares (CGLS) algorithm directly  to the problem $ b = \mA_\theta w$.  The CGLS algorithm determines a sequence of approximate solutions by solving the minimization problems in a nested sequence of Krylov subspaces,
\[
  w_k = {\rm argmin}\big\{ \|b - \mA_\theta w\|^2 \mid w \in {\mathscr K}_k(\mA_\theta b, \mA_\theta^\mT \mA_\theta)\big\},
\]  
where the Krylov subspace of order $k$ is defined as
\[
 {\mathscr K}_k(\mA_\theta b, \mA_\theta^\mT \mA_\theta) = {\rm span}\big\{ \mA_\theta b,(\mA_\theta^\mT \mA_\theta) \mA_\theta b,\ldots,
,(\mA_\theta^\mT \mA_\theta)^{k-1} \mA_\theta b\big\}.
\]
The CGLS method is a computationally efficient alternative to Tikhonov regularization, albeit not equivalent to it, provided that the iterations are terminated right before the amplified noise components start to overtake the solution. In the present problem, we stop the iterations as soon as one of the following conditions are satisfied,
\begin{equation}\label{stop CGLS1}
 \|b - \mA_\theta w_k\|^2 < m, 
 \end{equation}
 or 
 \begin{equation}\label{stop CGLS2}
 \|b - \mA_\theta w_{k+1}\|^2 + \|w_{k+1}\|^2 >   \|b - \mA_\theta w_{k}\|^2 + \|w_{k}\|^2 .
 \end{equation}
Condition (\ref{stop CGLS1}) is the classical Morozov discrepancy principle, recalling that for Gaussian white noise $E\sim{\mathcal N}(0,\mI_m)$, the noise level can be defined as 
\[
{\mathbb E} \|E\|^2  = {\rm trace}\big({\mathbb E} E E^\mT\big) =  {\rm trace}(\mI_m)= m.
\]
The condition (\ref{stop CGLS2}), on the other hand,  monitors the convergence of the original least squares problem whose solution the iterative process seeks to approximate. The updating procedure for $z$ is the same regardless of the hyperprior model for $\theta$, so it applies for both Phase I and Phase II of the hybrid algorithm. For further discussion of the iterative algorithm, we refer to \cite{CSbook}.

\subsection{Updating the metric}

Consider next the updating scheme of the metric. Assume that we have completed the IAS iteration, and calculated the pair $(z,\theta)$, where $z$ contains the coefficients of the representation of the approximate gradient $\nabla\widetilde u$ in the Whitney basis and $\theta$ is the vectors of the variances. More specifically, both the variance $\theta_j$ and the coefficient $z_j$  are associated to the edge $e_j$ of the mesh. Because of the discontinuity of the Whitney basis functions across the edges, pointwise evaluation of the approximate gradient is not well defined, therefore, to avoid ambiguities, we approximate the  gradient  over  the domain $\Omega$ by a continuous interpolant. For that purpose we chose the standard first order Cl\'{e}ment interpolant \cite{Clement}, which we briefly review here for the sake of completeness.

To interpolate the discontinuous function (\ref{nabla u}), we begin by considering its components separately, introducing the notation
\[
 \partial_1\widetilde u(x) =  \sum_{j=1}^{n_e} z_j  e_1\cdot  w_j = \sum_{j=1}^{n_e} z_j  w^1_j, \quad  \partial_2\widetilde u(x) =  \sum_{j=1}^{n_e} z_j e_2\cdot  w_j = \sum_{j=1}^{n_e} z_j  w^2_j,
\]
where $e_1$ $e_2$ are the canonical unit vectors.  Cl\'{e}ment interpolation proceeds by
\begin{enumerate}
\item finding a local patch approximation around each nodal point $v_i$ using orthogonal projections to local piecewise polynomials;
\item interpolating the values of the patch approximations at their central node as global piecewise polynomial.
\end{enumerate}

In the following, we describe the two steps of  first order Cl\'{e}ment interpolation in detail for $ \partial_1\widetilde u(x)$;  the procedure for other component is analogous. 

We define a  local patch $\omega_i$ where the approximation is carried out by picking a node $v_i$, and consider all triangles with $v_i$ as a vertex, see Figure~\ref{fig:patch Whitney}. We define
\[
 p_i(x) = \mbox{$L^2(\omega_i)$-orthogonal projection of $\partial_1\widetilde u|_{\omega_i}$ onto ${\mathbb P}_1(\omega_i)$,}
\] 
where ${\mathbb P}_1(\omega_i)$ is the space of piecewise first order polynomials over the patch  $\omega_i$. 
If $I_i$ is the index set of all nodes associated with the patch $\omega_i$, then 
\[
 p_i(x) = \sum_{j\in I_i} \alpha_j^i \psi_j(x).
\]
On the other hand,  the orthogonality of the projection is tantamount to the condition
\[
 \int_{\omega_i} \psi_\ell(x)\big( \partial_1\widetilde u (x) - p_i(x)\big) dx = 0, \quad \ell\in I_i.
\]
This set of conditions leads to a square system of linear equations that must be satisfied by the coefficients determining the interpolation polynomial in $\omega_i$.  More precisely,
\begin{eqnarray*}
 \int_{\omega_i} \psi_\ell(x)\big( \partial_1\widetilde u (x) - p_i(x)\big) dx &=& \sum_{j=1}^{n_e} z_j \int_{\omega_i} w_j^1(x)\psi_\ell(x)dx -\sum_{j\in I_i} \alpha^i_j \int_{\omega_i}\psi_j(x)\psi_\ell(x) dx \\
 &=&\sum_{j=1}^{n_e} b^i_{\ell j} z_j - \sum_{j\in I_i} g^i_{\ell j} \alpha^i_j  = 0,
 \end{eqnarray*}
which can be expressed in matrix form as
\[
 \mG^i \alpha^i = \mB^i z.
\]
This is the linear system that needs to be solved to find he coefficient vector $\alpha^i$.  

Once the first order polynomials $p_i$ have been found for each node, the Cl\'{e}ment interpolant is defined as
\[
 {\mathscr J}_{\rm Cle}\big(\partial_1\widetilde u\big)(x) = \sum_{i=1}^{n_v} p_i(v_i)\psi_i(x).
\] 
We denote the Cl\'{e}ment interpolation of the function $\nabla\widetilde u(x)$ by 
\[
\nabla^* u(x)=\big( {\mathscr J}_{\rm Cle}\big(\partial_1\widetilde u\big)(x), {\mathscr J}_{\rm Cle}\big(\partial_2\widetilde u\big)(x)\big).
\]

Having the interpolated function  $\nabla^* u(x)$ at our disposal, we can define the metric $\mG(x)$ at an arbitrary point $x\in\Omega$ as follows. Let $\delta>0$ be a fixed threshold value below which the interpolated gradient is considered negligible. Then:
\begin{enumerate}
\item If $|\nabla u^*(x)|>\delta$, the metric is anisotropic: We set
\[
 e_{\|}(x) = \frac{\nabla u^*(x)}{|\nabla u^*(x)|}, \quad e_{\perp}(x) \perp e_{\|}(x),
\]
and define
\begin{equation}\label{metric}
 \mG(x) = C|\nabla u^*(x)|^2 \left(e_{\|}(x)e_{\|}(x)^\mT + \frac 1{\alpha} e_\perp(x)e_\perp(x)^\mT\right),
\end{equation}
where $\alpha\geq 1$ is a given anisotropy ratio, and $C>0$ is a scaling factor that will be specified below.
\item If $|\nabla u^*(x)|\leq \delta$, the metric is isotropic, and we set
\begin{equation}\label{metric 2}
 \mG(x) = \delta \mI,
\end{equation}
where $\mI$ is the $2\times 2$ identity matrix. 
\end{enumerate} 
The metric defined in this manner produces relatively large and approximately equilateral elements when the interpolated gradient is small, while for large values of the interpolated gradient, the elements are smaller, with a prescribed eccentricity encoded in the anisotropy ratio $\alpha$.

In practice, the mesh generator of  FreeFEM \cite{Hecht} requires the specification of the approximate minimum and maximum  edge lengths. In order for the mesh generator to interface seamlessly with the metric defined here, we first specify $0<h_{\rm min}<h_{\rm max}$. Assuming that $h_{\rm min}$ represents the minimum edge length  in the direction parallel to $\nabla^* u$, we set
\[
 h_{\rm min}^2 = \frac 1{C M^2}, \quad M = {\rm max}\{ |\nabla^* u(v_j)|, 1\leq j\leq n_v\},
\]
from which we can compute the value of $C$. To select the threshold value $\delta$, we require that when $|\nabla^*u(x)| = \delta$, the edge length perpendicular to the gradient assumes the value $h_{\rm max}$, yielding 
\[
 h_{\rm max}^2 = \frac{\alpha}{C \delta^2},
\]
from which can be  compute the value of  $\delta$. Therefore, the algorithm requires the specification of only the three parameters $(h_{\rm min},h_{\rm max},\alpha)$, all having a direct geometric interpretation.

\section{Algorithm and computed examples}

We test the algorithm on two computed examples.  
For the sake of simplicity, we only consider linear inverse problems, and the examples presented here are similar to those discussed in the earlier paper \cite{CCPS} to allow a comparison of the results obtain with a Bayesian adaptive isotropic meshing algorithm proposed there with our novel anisotropic model. The problems are an X-ray fanbeam tomography with sparse illumination angle data, and an inverse source problem for Darcy flow.

\subsection{Tomography problem}

In this example, the data consists of $N$ illumination views of the target by a fan of rays as indicated in Figure~\ref{fig:Fanbeam}. The image area $\Omega$ is defined as the unit disc  centered at the origin, and we assume that each illumination comprises $n$ uniformly spaced rays traversing the image area. Thus, the dimension of the data is $m = Nn$. In the computed example, we set $N = 15$ and $n = 300$, so the data have dimension $m = 45\,000$.
\begin{figure}
\centerline{\includegraphics[width=7cm]{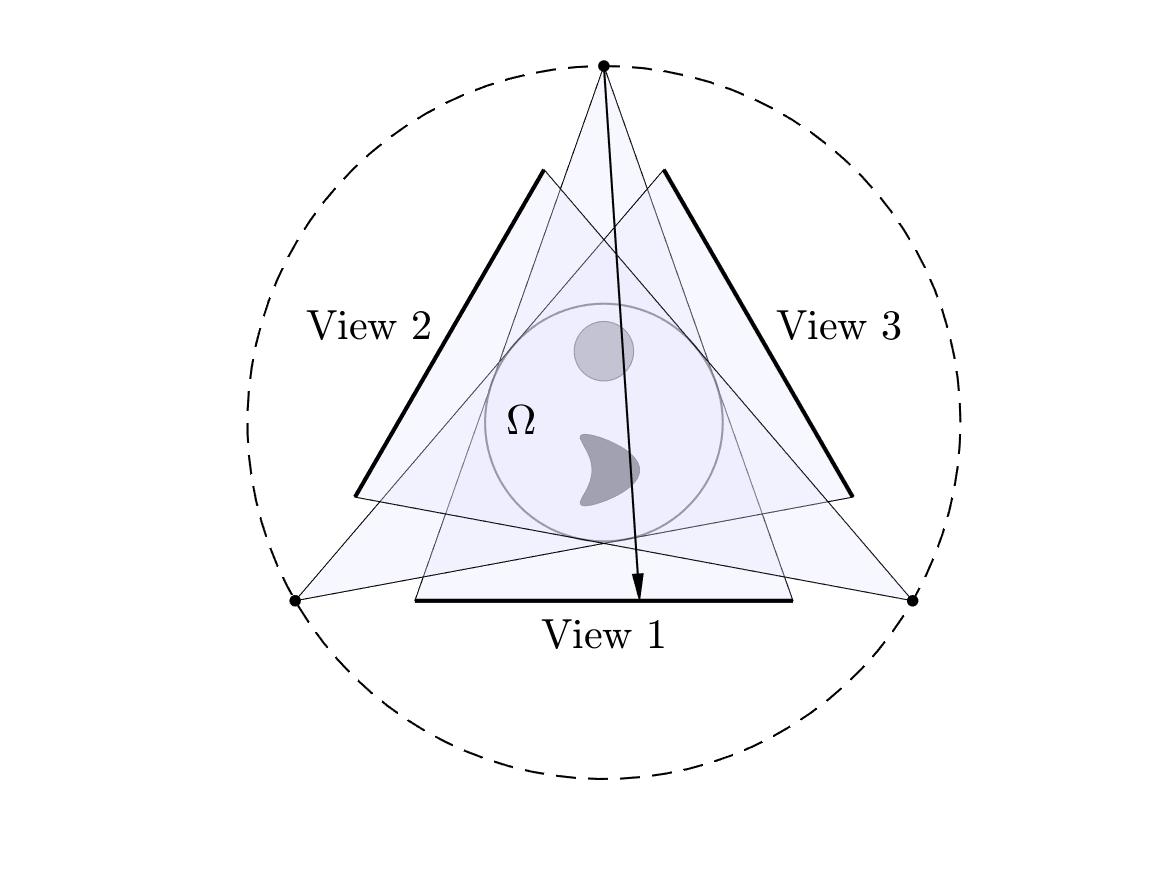}}
\caption{\label{fig:Fanbeam} A cartoonish rendition of the geometric setup for the first computed example.  Here, for simplicity, the number of views is limited to three. In the computed example, 15 equally spaced illumination directions are used. The number of traversing rays per view is 300, leading to data with dimension $m=45\,000$.}
\end{figure}

The attenuation of the intensity $I$ of the radiation along an infinitesimal line segment $ds$ is given by the Beer-Lambert law \cite{Natterer},
\[
 d I = - \mu ds,
\]  
where $\mu = \mu(x)$ is the absorption coefficient.
Thus, if a ray is represented in a form  $x = x(s)\in\Omega$ parametrized by its arc length $s$, the attenuated intensity at the end of a ray of length $L$ is obtained by the formula
\[
 I_L = I_0 \, {\rm exp}\left( - \int_0^L \mu(x(s)) ds\right).
\]
We define the noiseless data $b^*\in \R^{Nn}$ componentwise as
\[
 b_k^* =  - \log  \frac{I_{L_k}}{I_0}  =  \int_0^{L _k} \mu(x_k(s)) ds, \quad 0\leq s\leq L_k,
\] 
where $x_k(s)$ is the parametrization of the $k$th ray. We assume for simplicity  that the noise in the components $b_k^*$ can be approximated by additive Gaussian scaled white noise,
\[
 b = b^* + w, \quad w \sim{\mathcal N}(0, \sigma^2 \mI_{m}),\quad m = Nn.
\] 
In the numerical simulation,  we assume that the target is a translucent body with two absorbing inclusions with different, yet constant absorption coefficients. 
 
To highlight different features of the algorithm, we run three protocols. In the first run, the noiseless data are corrupted by additive Gaussian scaled white noise, with standard deviation $\sigma$ equal to 4\% of the maximum of the noiseless signal, which corresponds to $\sigma = 0.0472$.  In the first phase of the hybrid IAS algorithm, which uses the gamma hyperprior  parameter, we set  $\eta = 0.001$  and $\vartheta^*_1 = 0.05$.  In this example, the sensitivity scaling of the vector $\vartheta$ is not necessary, thus, for simplicity, the components of the vector $\vartheta$ are all set equal.  When the relative change  of the vector $\theta$ is below a threshold value of 5\%, the iterations in Phase I are stopped and as we move to Phase II,  the hypermodel  changes to a generalized gamma hyperprior with $r=1/2$. The hyperparameter values for Phase II are set automatically by using the matching conditions. Similarly to Phase I, the IAS iterations are stopped when the relative change of $\theta$ is below 5\%. However,  in both phases, we set the maximum number of IAS iterations to 15. In the remeshing algorithm, the mesh size parameters are  $h_{\rm min} = 0.01$ and $h_{\rm max} = 0.1$, while the anisotropy factor is $\alpha = 12$.

We start the iterations by defining an isotropic mesh with uniform metric, with mesh size parameter $h_{\rm init} = 0.05$. The initial mesh is shown in Figure~\ref{fig:high noise tomo} top row, left panel, superimposed with the true target. The absorption coefficient in the circular inclusion of the true target is $\mu = 1.2$, and in the kite-shaped inclusion $\mu = 1.0$, while the background is for simplicity assumed to have a vanishing absorption coefficient. 

We run four full iterations of the IAS-remeshing algorithm, and show the results at the end of each iteration in Figure~\ref{fig:high noise tomo}. The left panel of each row shows the current tessellation, and the middle and right panels show the corresponding reconstructions. Numerical experiments show that after four iterations, the reconstructions and the meshes do not change significantly.
\begin{figure}[ht!]
\centerline{\includegraphics[width=13cm]{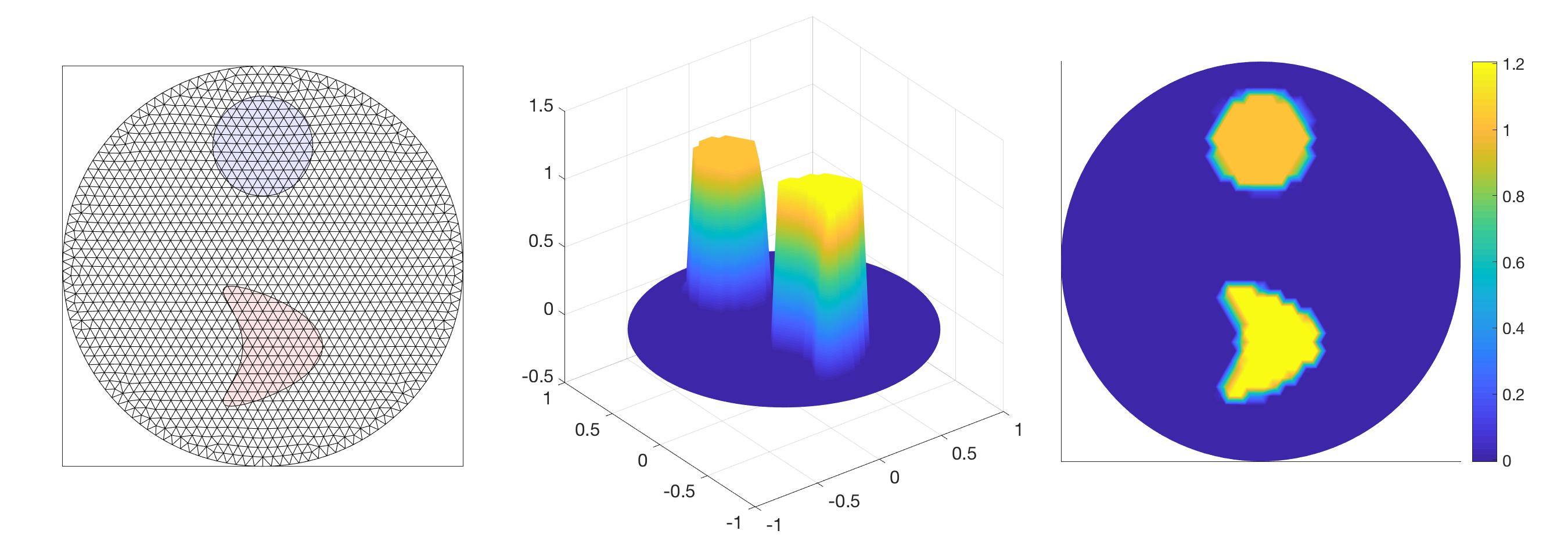}}
\centerline{\includegraphics[width=13cm]{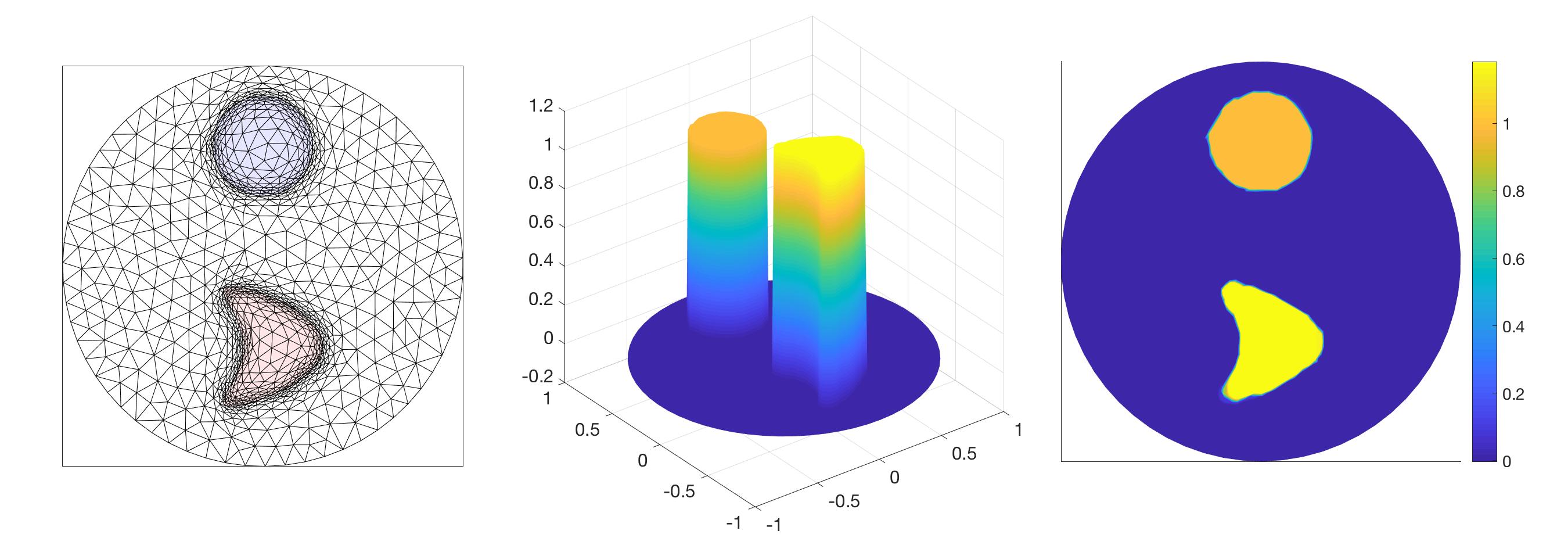}}
\centerline{\includegraphics[width=13cm]{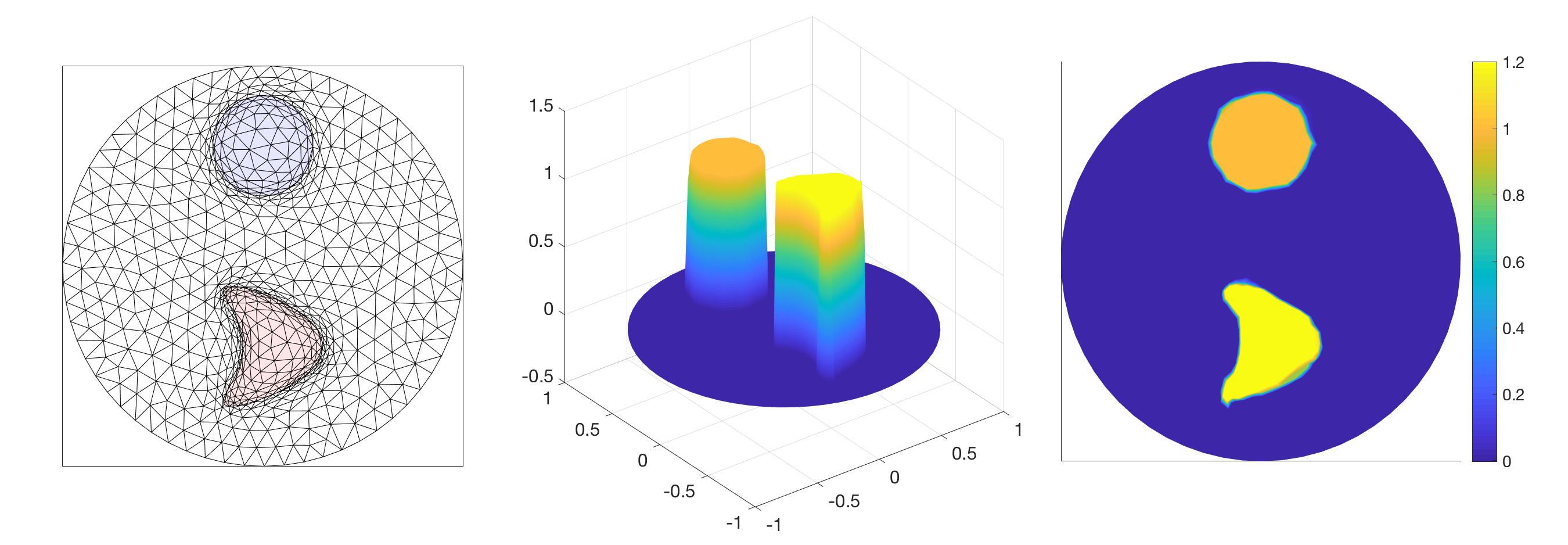}}
\centerline{\includegraphics[width=13cm]{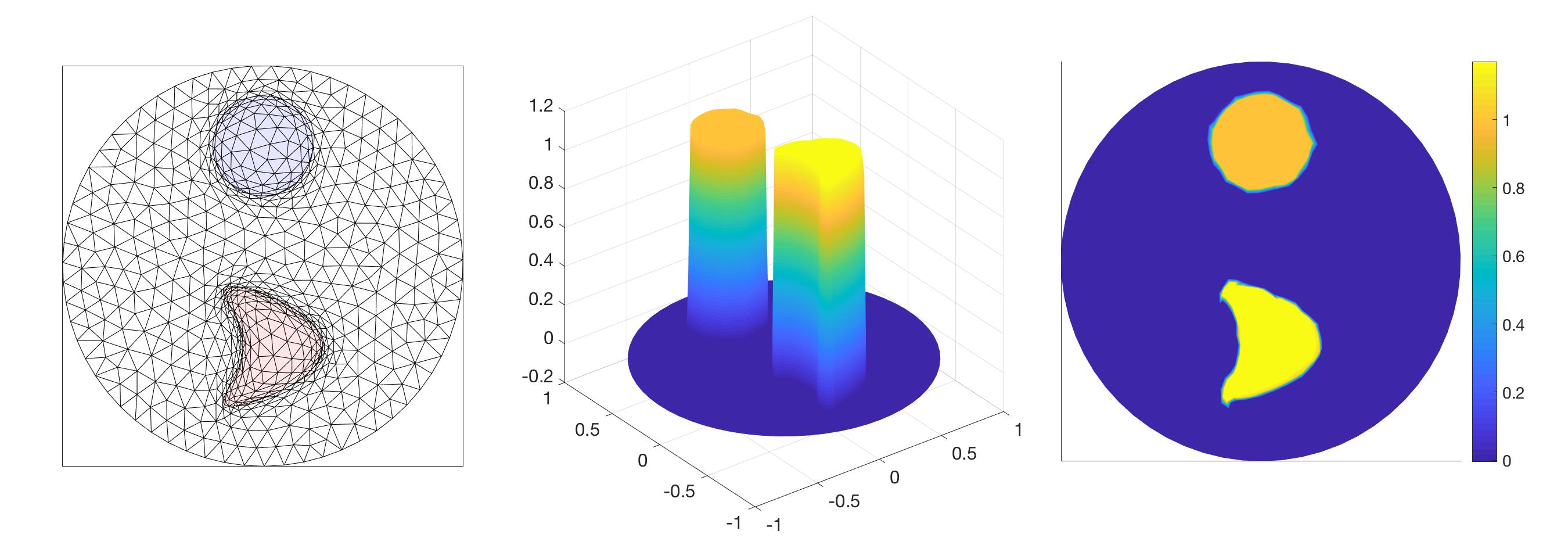}}
\caption{\label{fig:high noise tomo} Four iterations of the hybrid IAS-remeshing algorithm. Here, the data are contaminated by high noise level  of 4\% of the maximum of the noiseless signal.}
\end{figure}

To assess the efficiency of the algorithm, we consider first the computational complexity of the IAS algorithm, and in particular, the number of inner and outer iterations.  Table~\ref{tab:CGLS ex1} lists the number of outer IAS iterations of both phases if the hybrid IAS algorithm, as well as the number of the CGLS iterations. We recall that the CGLS iterations are stopped when one of the stopping criteria (\ref{stop CGLS1}) or (\ref{stop CGLS2}) is satisfied.
 We observe that each IAS update requires fewer than 20 CGLS iterations, making the IAS algorithm very fast. 
\begin{table}[ht!] 
\centerline{
\begin{footnotesize}
\begin{tabular}{c|cccccccccc|ccccc}
\multicolumn{1}{c|}{} & \multicolumn{10}{c|}{{\bf Phase I}} & \multicolumn{5}{c}{{\bf Phase II}} \\
&1&2&3&4&5&6&7&8&9&10&1&2&3-8 &9-11&12-15\\
\hline
I   &18  & 14 & 10  & 9 & 8 & 8 & 8   &     &    &    &           4 & 3 & 3 &   &       \\
II  & 18 & 14 & 12  & 9   & 9 & 8 & 7 & 7  & 7 &     &          4 & 4 & 3 & 3 &  3   \\ 
III & 18 &  14 & 12 & 10 & 9 & 9 & 8 &     &    &     &          5 & 4 & 3 & 3 &      \\ 
IV & 18 &  14 & 11 & 9  &  9 &9 &  8  & 8 & 7 & 7  &          4 & 3 & 3 & 3 &  3  \\ 
 \end{tabular}
\end{footnotesize}
}
\caption{\label{tab:CGLS ex1} Number of CGLS iterations in the two phases corresponding to the first simulation. The rows refer to the full iterations consisting of the full hybrid IAS updates and remeshing, and the columns correspond to IAS iterations with gamma hyperprior (Phase I) and generalized gamma hyperprior with $r = 1/2$ (Phase II). The maximum number of IAS iterations allowed was 15, and a missing number indicates that the stopping criterion of relative change of $\theta$ dropping under 5\% was satisfied before the maximum number of iterations was reached.}
\end{table}

Before discussing the complexity of the algorithm measured in terms of mesh sizes and computing times, we address first the question whether the  mesh refinement is fine enough to ignore the modeling error. Recall that the artifacts due to the modeling error become particularly significant when the data are of high quality, in which case the discretization error dominates the exogenous noise. Therefore, we test the algorithm using the same parameter values for the IAS and the mesh refinement as in the previous run, but decrease the noise level in the data by setting the standard deviation of the noise at $1\%$  of the maximum of the noiseless signal, or equivalently, setting $\sigma = 0.0118$. To avoid that the first IAS iteration ends up fitting the data to the discretization error, in the first iteration we artificially inflate the standard deviation in the likelihood and set $\sigma = 0.3\times h_{\rm init} = 0.015$. The motivation for this choice is that in \cite{CCPS}, it was found through numerical simulations that in the tomography problem with piecewise constant inclusion, the modeling error level increases roughly linearly as a function of the mesh parameter with the empirical coefficient 0.3. No artificial inflation is applied after the first iteration.   

\begin{figure}[ht!]
\centerline{\includegraphics[width=13cm]{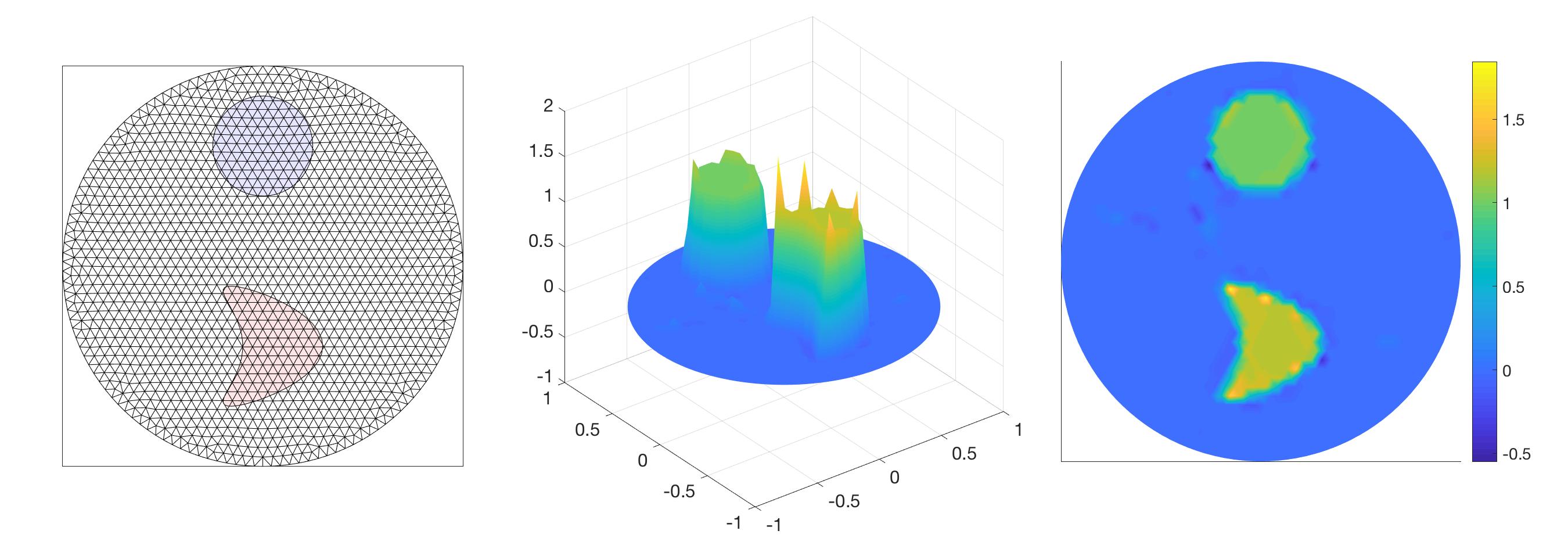}}
\centerline{\includegraphics[width=13cm]{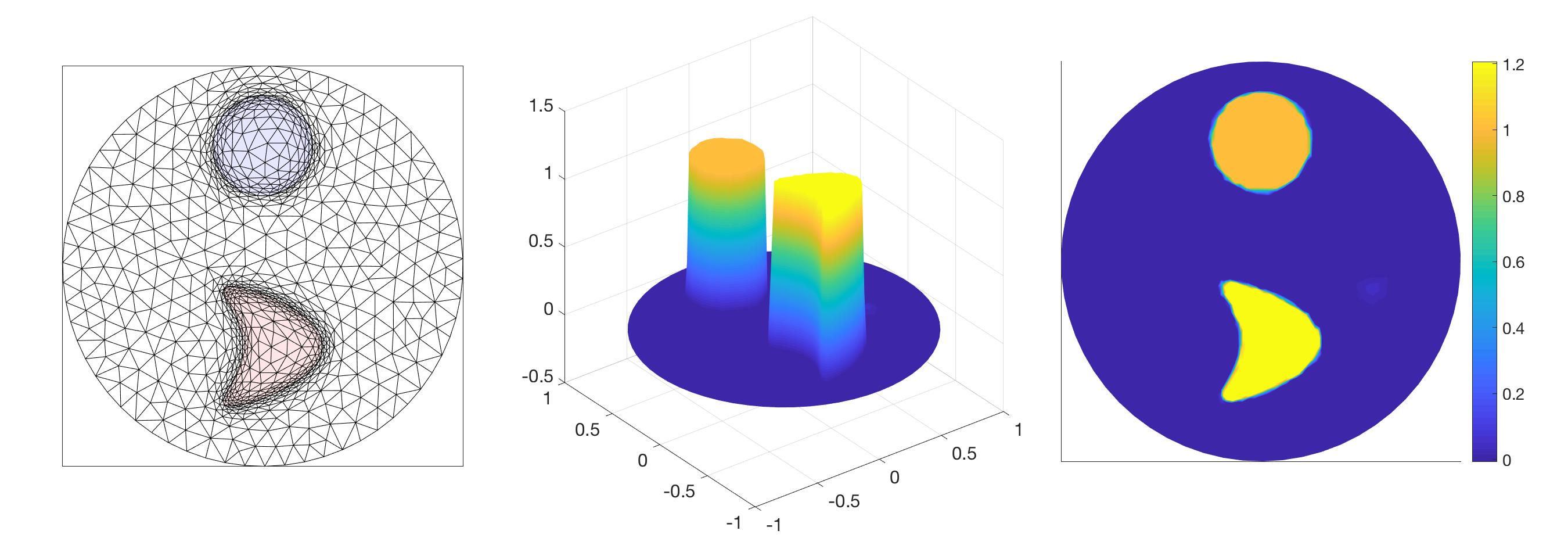}}
\centerline{\includegraphics[width=13cm]{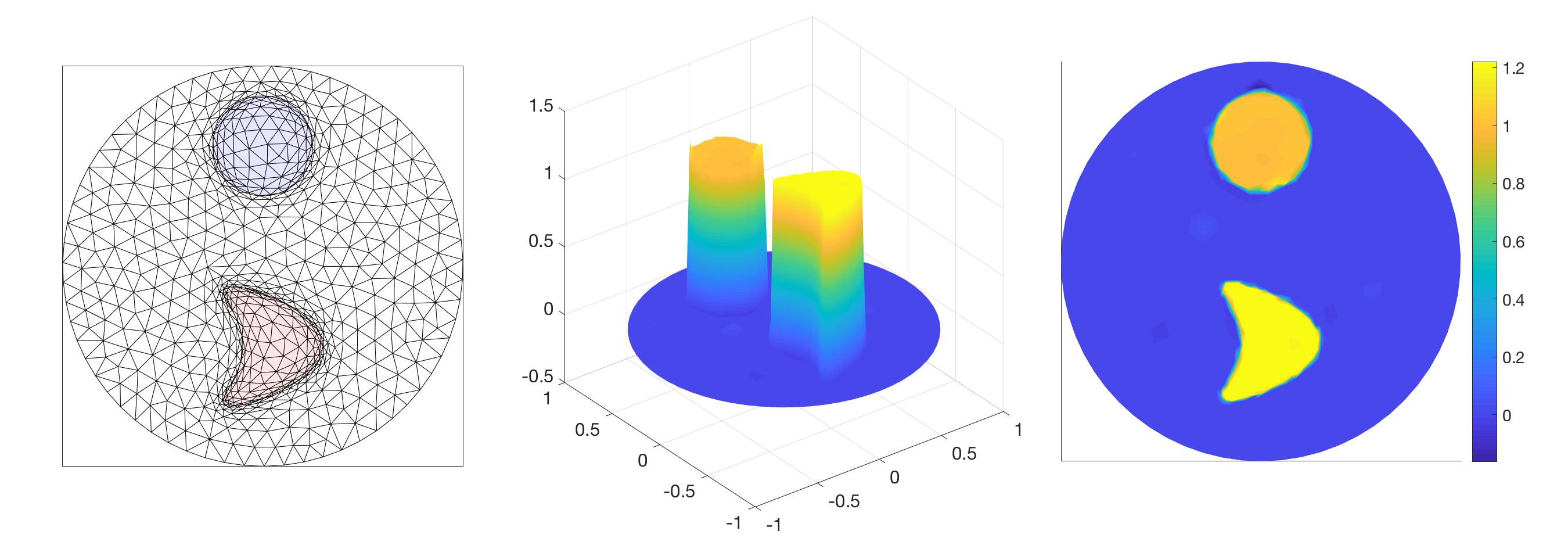}}
\centerline{\includegraphics[width=13cm]{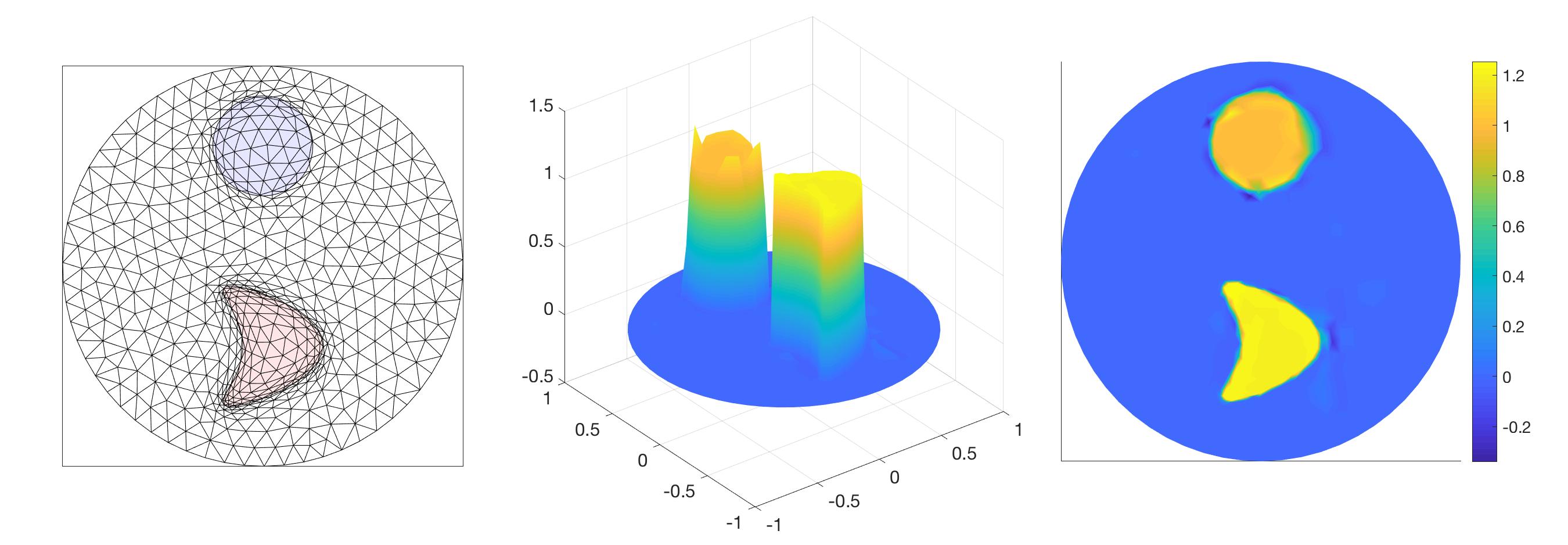}}
\caption{\label{fig:low noise tomo} Four iterations of the hybrid IAS-remeshing algorithm with low noise of 1\% of the maximum of the noiseless signal. In this case, the mesh refinement is the same as in the high noise example, making the approximation error the dominant part. The effect of the modeling error is visible in the reconstructions.}
\end{figure}

The results of the first iteration rounds are shown in the top row of Figure~\ref{fig:low noise tomo}. In the first reconstructions, the  amplified modeling errors are clearly visible, and even in fourth iteration boundary artifacts persist, yielding worse reconstructions than those obtained with higher noise level. This indicates that the refinement parameter $h_{\rm min}$ in the metric was not chosen small enough.
The number of the CGLS iterations within the IAS iterations are listed in Table~\ref{tab:CGLS ex2}. As expected, the lower noise level requires a better fit to the data, thus more inner iterations are required.
\begin{table} 
\centerline{
\begin{footnotesize}
\begin{tabular}{c|ccccccc|cccccccc}
\multicolumn{1}{c|}{} & \multicolumn{7}{c|}{{\bf Phase I}} & \multicolumn{8}{c}{{\bf Phase II}} \\
&1&2&3&4&5&6&7&1&2&3&4&5&6&7&8 \\
\hline
I  &44  & 31 & 28 & 27 & 26 & 24 & 24& 11 & 9 & 10 & 9 &    &    &    &    \\
II & 59 & 34 & 29 & 28 & 26 &26 &      & 10 & 7 &   6 & 6 &  6  & 6 & 6 & 6 \\ 
III & 60 & 37 & 31 & 29 & 29 &  &      & 10 & 8 &   7 & 7 &  8  & 7 &   &   \\ 
IV & 61 & 38 & 35 & 35 & 35 &34 & 33     & 14 & 12 &  11 & 11 &  11  &  &  &  \\ 
 \end{tabular}
\end{footnotesize}
}
\caption{\label{tab:CGLS ex2} Number of CGLS iterations in the two phases corresponding to the second simulation. The rows correspond to the full iterations consisting of the full hybrid IAS updates, and the columns correspond to IAS itrations with gamma hyperprior (Phase I) and generalized gamma hyperprior with $r = 1/2$ (Phase II). The maximum number of IAS iterations allowed was 15, and a missing number indicates that the stopping criterion was reached before this limit.}
\end{table}

The third run is almost identical to the second one, except for the halving of the mesh size parameter in the remeshing, which is set to  $h_{\rm min} = 0.005$, that is, the minimum mesh edge is halved. The results of the four iterations of the algorithm are shown in in Figure~\ref{fig:low noise tomo fine}.  This time, there are no boundary artifacts around the inclusions, and the algorithm approximates well the inclusions. The number of CGLS iterations is given in tabular form in Table~\ref{tab:CGLS ex3}.

\begin{figure}[ht!]
\centerline{\includegraphics[width=13cm]{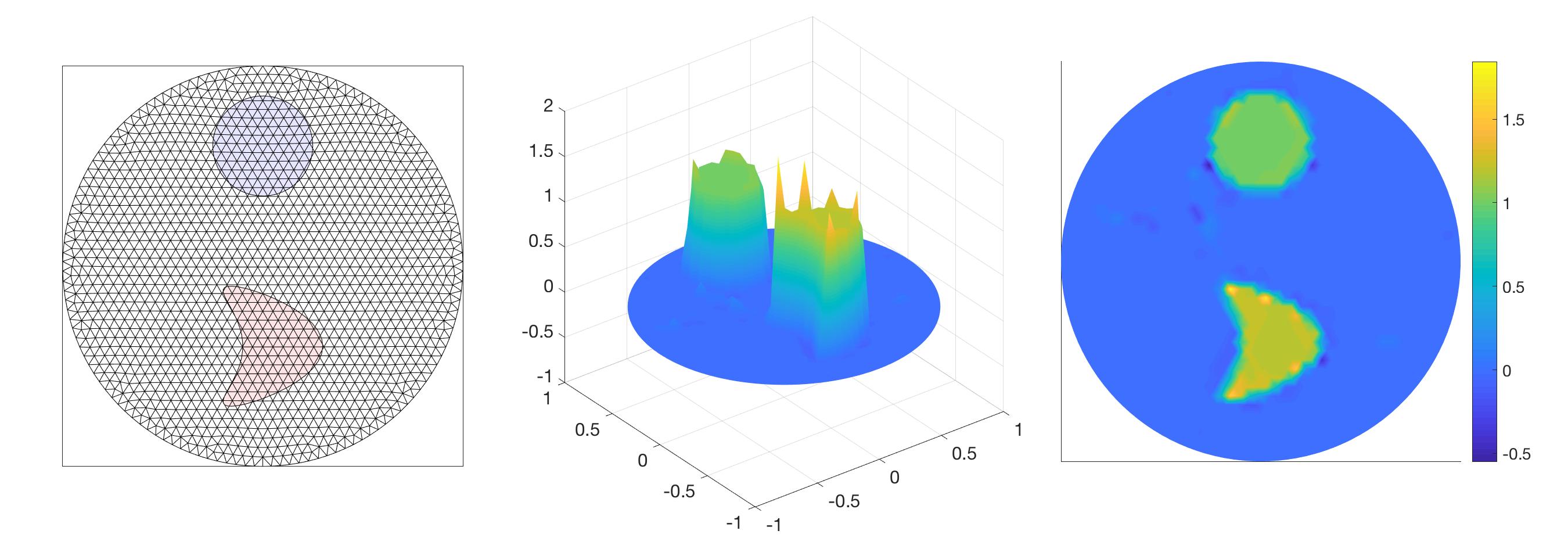}}
\centerline{\includegraphics[width=13cm]{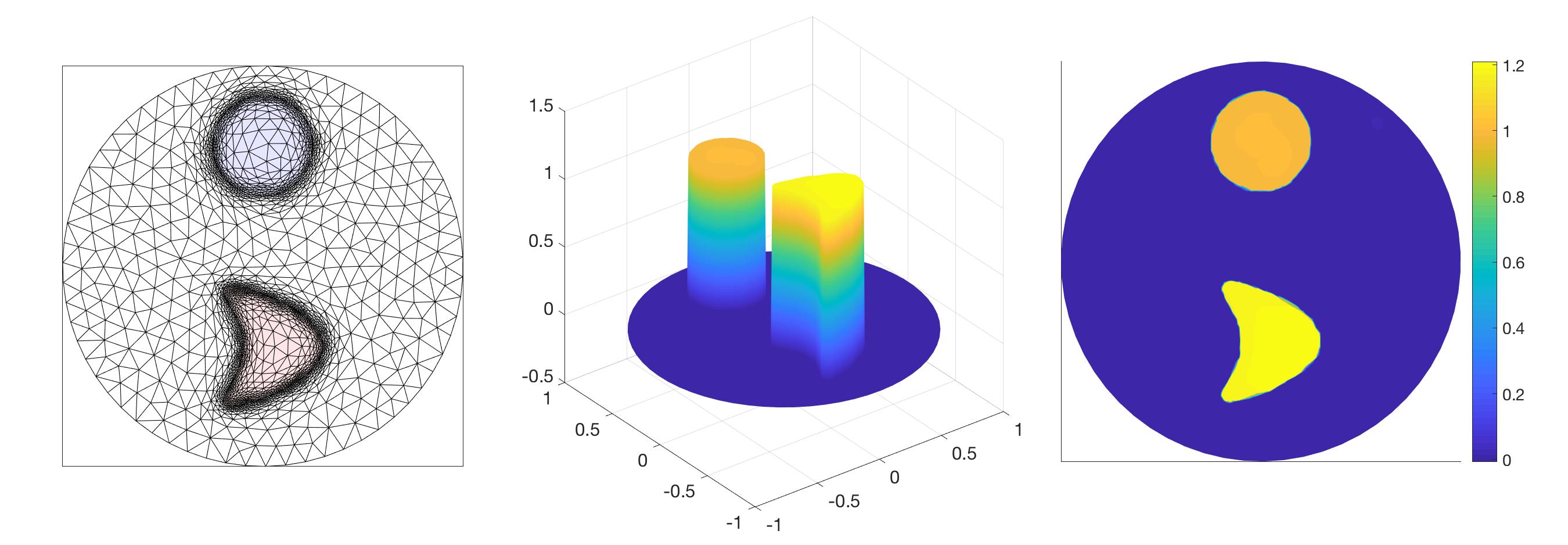}}
\centerline{\includegraphics[width=13cm]{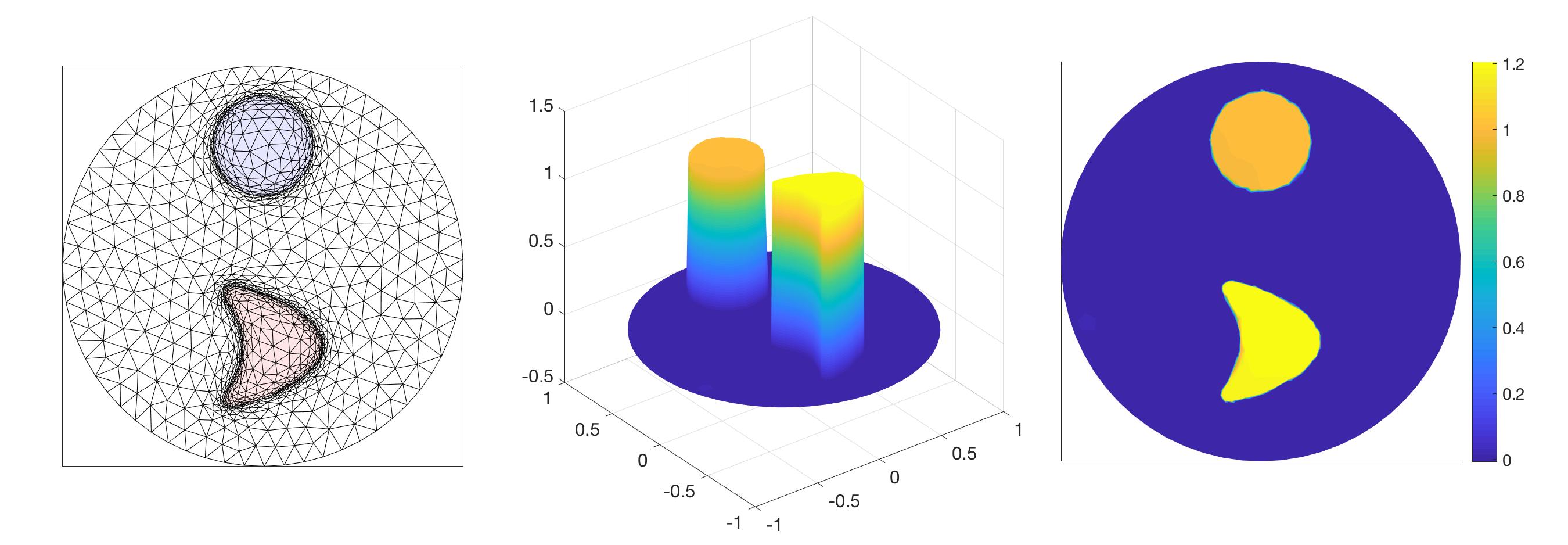}}
\centerline{\includegraphics[width=13cm]{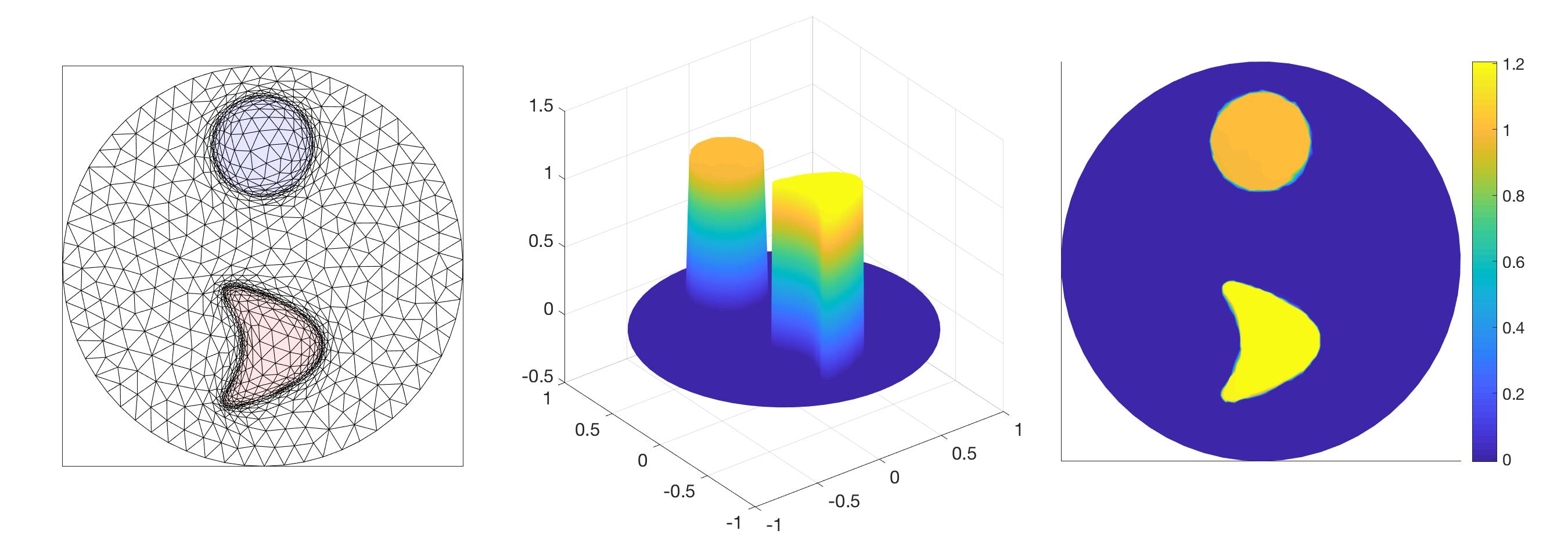}}
\caption{\label{fig:low noise tomo fine} Four iterations of the hybrid IAS-remeshing algorithm with low noise of 1\% of the maximum of the noiseless signal. Here, the minimum edge length in the mesh refinement is one half of the previous example, reducing the modeling error below the additive noise level, yielding a better reconstruction. }
\end{figure}

\begin{table} 
\centerline{
\begin{footnotesize}
\begin{tabular}{c|ccccccccc|ccccccccccccccc}
\multicolumn{1}{c|}{} & \multicolumn{9}{c|}{{\bf Phase I}} & \multicolumn{7}{c}{{\bf Phase II}} \\
&1&2&3&4&5&6&7&8&9&1&2&3&4&5&6-14&15 \\
\hline
I  &44  & 31 & 28 & 27 & 26 & 24 & 24&      &       & 11 & 9 & 10 & 9 &    &    &     \\
II & 59 & 28 & 21 & 18 & 17 &16 & 16 & 15 & 15  &   8 & 7 &  7  & 7 & 7 & 6 & 6 \\ 
III & 58 & 32 & 22 & 18 & 17 &16 &  &  &   &   9 & 7 &  6  & 6 & 6 & 6 & 6      \\   
IV & 59 & 32 & 22 & 19 & 19 &18 & 18 &  &   &   9 & 7 &  7  & 7 & 7 & 6 &  \\  
 \end{tabular}
\end{footnotesize}
}
\caption{\label{tab:CGLS ex3} Number of CGLS iterations in the two phases corresponding to the third simulation. The rows correspond to the full iterations consisting of the full hybrid IAS updates, and the columns correspond to IAS itrations with gamma hyperprior (Phase I) and generalized gamma hyperprior with $r = 1/2$ (Phase II). The maximum number of IAS iterations allowed was 15, and a missing number indicates that the stopping criterion was reached before this limit.}
\end{table}

To have closer look at the effect of the modeling error, in Figure~\ref{fig:profiles} we plot the estimated absorption coefficients corresponding to experiments 2 and 3 along the vertical diagonal through the disc $\Omega$. The first reconstructions perfectly coincide, as all parameters are the same. The second reconstructions are relatively good in both examples, however, in the third and fourth reconstruction, when the algorithm starts to sparsify the mesh where the gradient is believed to be small, boundary artifacts with the coarser mesh appear. These artifacts can be attributed to the discretization error. 

\begin{figure}
\centerline{
\includegraphics[width=15cm]{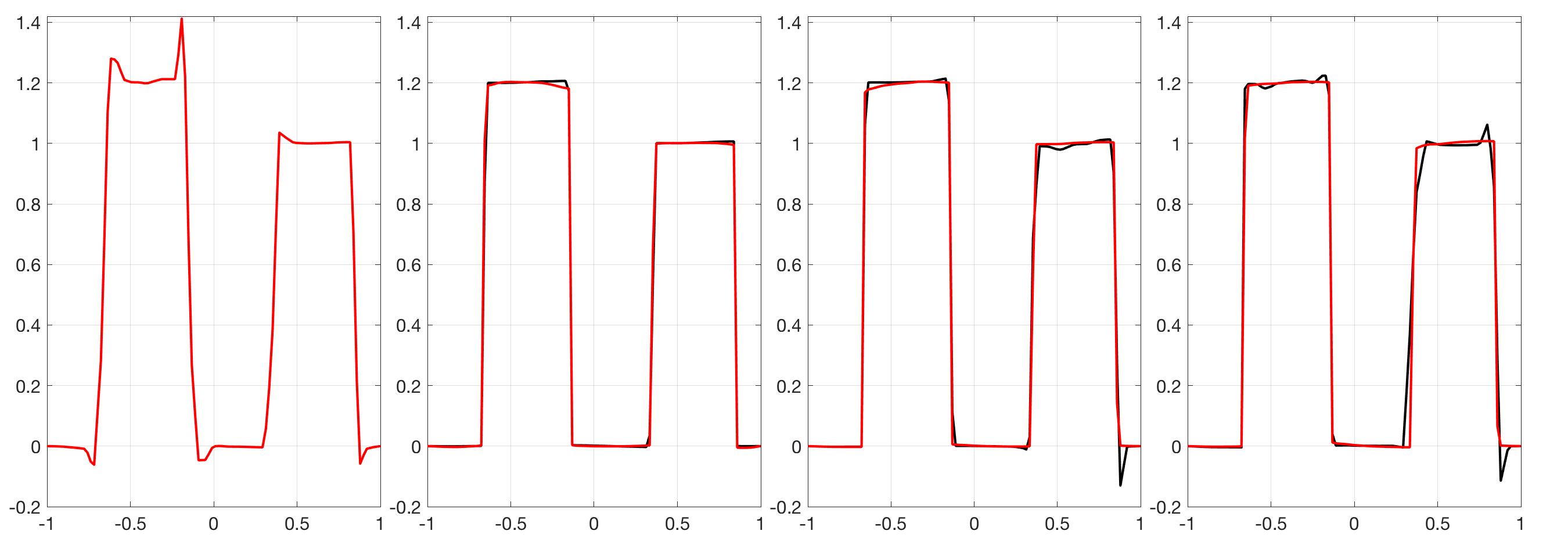}
}
\caption{\label{fig:profiles} Profiles of the reconstructions along the vertical diagonal of the disc $\Omega$. The black curve corresponds to the coarser discretization (Experiment 2), the red one to the finer discretization (Experiment 3). The four plots correspond to the four full iterations of the IAS and remeshing algorithm.}
\end{figure}

Finally, we address the question of complexity of the algorithm in terms of number of elements and computing times. Table~\ref{tab:times} shows the mesh sizes of each full iteration (I--IV) of the algorithm, together with the computing times. We observe that the first remeshing increases the number of elements, creating an overly dense mesh around a wide area around the object boundaries, after which the algorithm learns the location of the discontinuities and reduces the complexity to only the necessary mesh.

\begin{table}
\centerline{
\begin{footnotesize}
\begin{tabular}{c|cc|cc|cc}
\multicolumn{1}{c|}{} &\multicolumn{2}{c|}{{\bf Experiment 1}} &\multicolumn{2}{c|}{{\bf Experiment 2}}& \multicolumn{2}{c}{{\bf Experiment 3}}\\
 & $n_t$ & time $[{\rm s}]$ &  $n_t$ & time $[{\rm s}]$ &  $n_t$ & time $[{\rm s}]$ \\
\hline
I  & 2774 & 22.5 & 2774 & 20.4  & 2774 & 19.7 \\
II &1931 & 16.2 & 1691  & 10.5  & 3901 & 67.9 \\
III &1285 & 6.2 & 1251   & 5.1   & 2019 & 12.1 \\
IV &1223 & 7.5 & 1127  & 5.6    & 1763 & 14.4 \\
 \end{tabular}
\end{footnotesize}
}
\caption{\label{tab:times} The number of elements in the meshes, and the total computing time of the full updating of the estimates and the meshes. The computations are performed with Matlab on a standard laptop (MacBookPro 3.1 GHz intel Core i5). The recomputing of the forward model is not included in the times.}
\end{table}

For comparison, we run the three experiments setting the anisotropy parameter to $\alpha=1$, that is, the refined meshes are generated by using isotropic metric as in \cite{CCPS}. Table~\ref{tab:times iso} shows the mesh sizes and computing times, underlining the advantage of allowing anisotropy in the meshes. In particular, the first remeshing that is always the densest one creates a slowdown of the algorithm by a factor of 4 to 8. Finally, we point out that if the original mesh parameter $h=0.05$ would be refined homogeneously over the entire computational domain to yield a mesh parameter $h=h_{\rm min} = 0.01$ as in experiments 1 and 2, each element would be split into 25 subelements, yielding a mesh with $n_t = 69\,350$, resulting in significant computational complications.

\begin{table}
\centerline{
\begin{footnotesize}
\begin{tabular}{c|cc|cc|cc}
\multicolumn{1}{c|}{} &\multicolumn{2}{c|}{{\bf Experiment 1}} &\multicolumn{2}{c|}{{\bf Experiment 2}}& \multicolumn{2}{c}{{\bf Experiment 3}}\\
 & $n_t$ & time $[{\rm s}]$ &  $n_t$ & time $[{\rm s}]$ &  $n_t$ & time $[{\rm s}]$ \\
\hline
I  & 2774 & 22.2 & 2774 & 19.8  & 2774 & 20.3 \\
II &3853 & 68.6 & 3229  & 42.8  & 9863 & 585.2 \\
III &1581 & 10.4 & 1765   & 10.1   & 3727 & 61.0 \\
IV &1661 & 10.0 & 1705  & 9.9    & 2803 & 33.5 \\
 \end{tabular}
\end{footnotesize}
}
\caption{\label{tab:times iso} The number of elements in the meshes, and the total computing time of the full updating of the estimates and the meshes using isotropic metric, $\alpha = 1$. The recomputing of the forward model is not included in the times.}
\end{table}

\subsection{Inverse source problem}

Let $\Omega = [0,1]\times[0,1]\in\R^2$ denote the domain of interest, $u\in L^\infty(\Omega)$ a source function, and let $f\in H^1(\Omega)$ be the weak solution of the Dirichlet problem,
\[
 -\Delta f = u, 
\]
with a mixed boundary condition: Defining the vertical and horizontal boundary parts as
\begin{eqnarray*}
 \partial\Omega_v &=& \big\{x = (0,x^{(2)})\mid 0<x^{(2)}<1\big\}\cup \big\{x = (1,x^{(2)})\mid 0<x^{(2)}<1\big\}, \\
  \partial\Omega_h &=& \big\{x = (x^{(1)},0)\mid 0<x^{(1)}<1\big\}\cup \big\{x = (x^{(1)},1)\mid 0<x^{(1)}<1\big\},
\end{eqnarray*}
we impose
\[
 f\big|_{\partial\Omega_v} = 0, \quad \frac{\partial f}{\partial n}\bigg|_{\partial\Omega_h} = 0.
\] 
In this example the data consists of noisy observations of the solution $f$ at the nodes of a regular grid of size $20\times 20$, 
\[
 b_{j} = f(y_{j}) + \varepsilon_{j}, \quad y_{j} = \left(\frac{i-1/2}{20},\frac{k -1/2}{20}\right)\in\Omega, \quad j = 20(i-1) + k,\quad 1\leq i,k\leq 20,
\] 
The inverse problem is to estimate the source term $u$ from this data. As in the previous example, we assume that the source term is known to satisfy  $u\big|_{\partial \Omega} = 0$.

Given a tessellation ${\mathscr T}_h$,  we  approximate the source term $u$ in terms of the first order Lagrange basis functions  $\{\psi_j\}$ as in (\ref{u discr}). To approximate $f$, we use the second order piecewise polynomial Lagrange basis functions  $\{\varphi_j\}$  associated to nodes  obtained by augmenting the vertices $\{v_j\}$ by the midpoints of the edges of the triangles. Denoting by $I_{\rm free}$ the indices to nodes  satisfying $0<v^{(1)}_j<1$, the Galerkin approximation of the weak form leads to the linear system for the nodal values $\{f_j\}_{j\in I_{\rm free}}$ of the hydraulic potential, arranged in a vector $f_h$,
\[
 \mG_h f_h = \mM_h u_h,
\]
where the entries of the corresponding matrices  are given by
\[
 (\mG_h)_{jk} = \int_\Omega \nabla \varphi_j \cdot \nabla \varphi_k dx, \quad  j,k\in I_{\rm free}, 
\]
\[ 
 \quad  (\mM_h)_{j\ell} = \int_\Omega   \varphi_j  \psi_\ell dx, \quad j\in I_{\rm free}, \quad 1\leq \ell\leq n_v.
\]
The matrix $\mG_h$ is symmetric positive definite, thus the vector $f_h$ can be written as
\[
 f_h = \mG_h^{-1}\mM_h u_h.
\] 
Assuming that the data consist of evaluations of $f$ at points $p_i\in\Omega$, $1\leq i\leq m$, we have the discretized model for the observation,
\[
 b = \mA_h u_h,
\]
where the matrix $\mA_h$ is given by
\[
 \mA_h = \mP_h    \mG_h^{-1}\mM_h, \quad (\mP_h)_{jk} = \varphi_k(y_j).
\] 
 
We generate the data using a fine hexagonal mesh with mesh size parameter $h_{\rm fine} = 0.02$. The data are corrupted by scaled Gaussian white noise with standard deviation $\sigma = 3\times 10^{-4}$, corresponding to a noise level of  $10^{-4}$ of the maximum of the noiseless data.

For the inverse problem,  we run the algorithm starting with a regular hexagonal coarse mesh with mesh size parameter $h = 0.05$, performing four full iterations of IAS and subsequent remeshing.  As in the previous example, we run the hybrid IAS algorithm with $r_1=1$ and $r_2 = 1/2$.   Unlike in the tomography problem, here it is important for the quality of the reconstruction to use the sensitivity scaling (\ref{sens scaling}), due to the fact that the mesh edges are at different distances from the data points, and without the scaling adjustment, the coefficients $z_j$ related to edges near the data points would be favored. The hyperparameters for the IAS algorithm were set to $\eta = \beta_1 - 3/2 = 0.001$ and $\vartheta_1^* = 0.05$, and the remeshing parameters were chosen as $h_{\rm min} = 0.003$, $h_{\rm max} = 0.05$, and the anisotropy factor is $\alpha =12$.
Again, in the first iteration round, we inflate the variance of the likelihood to compensate for the unknown discretization error, using the same inflation factor as in the previous example, $\sigma = 0.3\times h$.

 \begin{figure}
\centerline{\includegraphics[width=13cm]{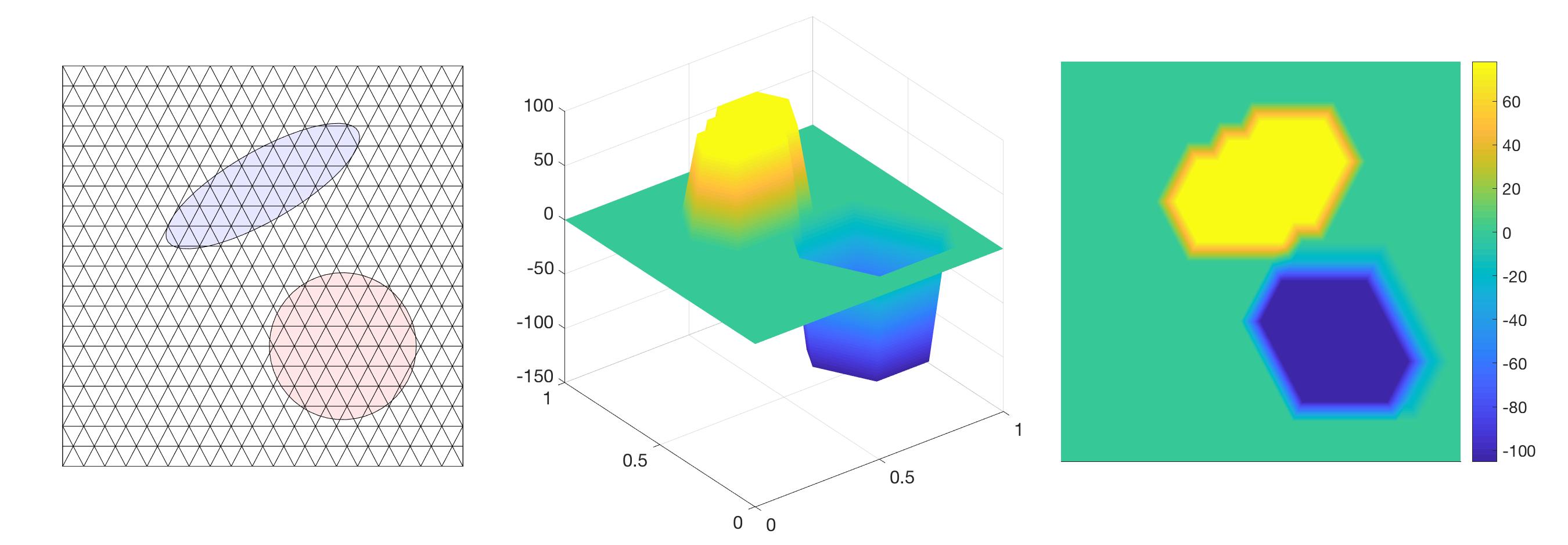}}
\centerline{\includegraphics[width=13cm]{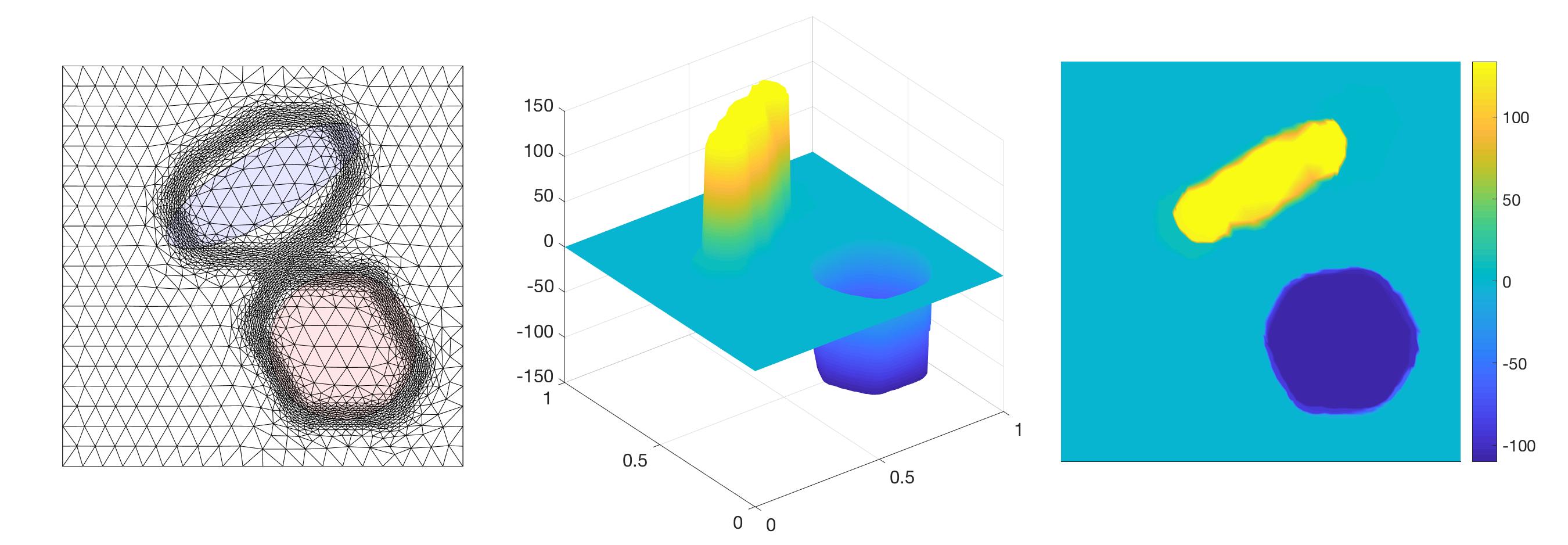}}
\centerline{\includegraphics[width=13cm]{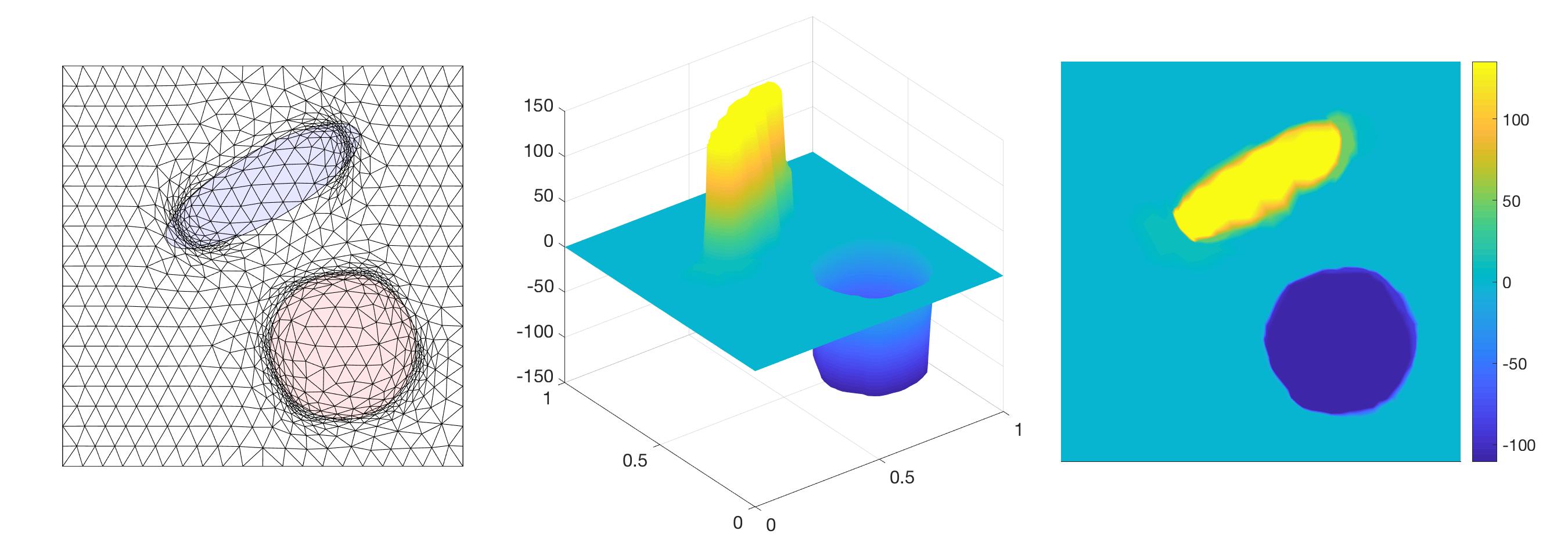}}
\centerline{\includegraphics[width=13cm]{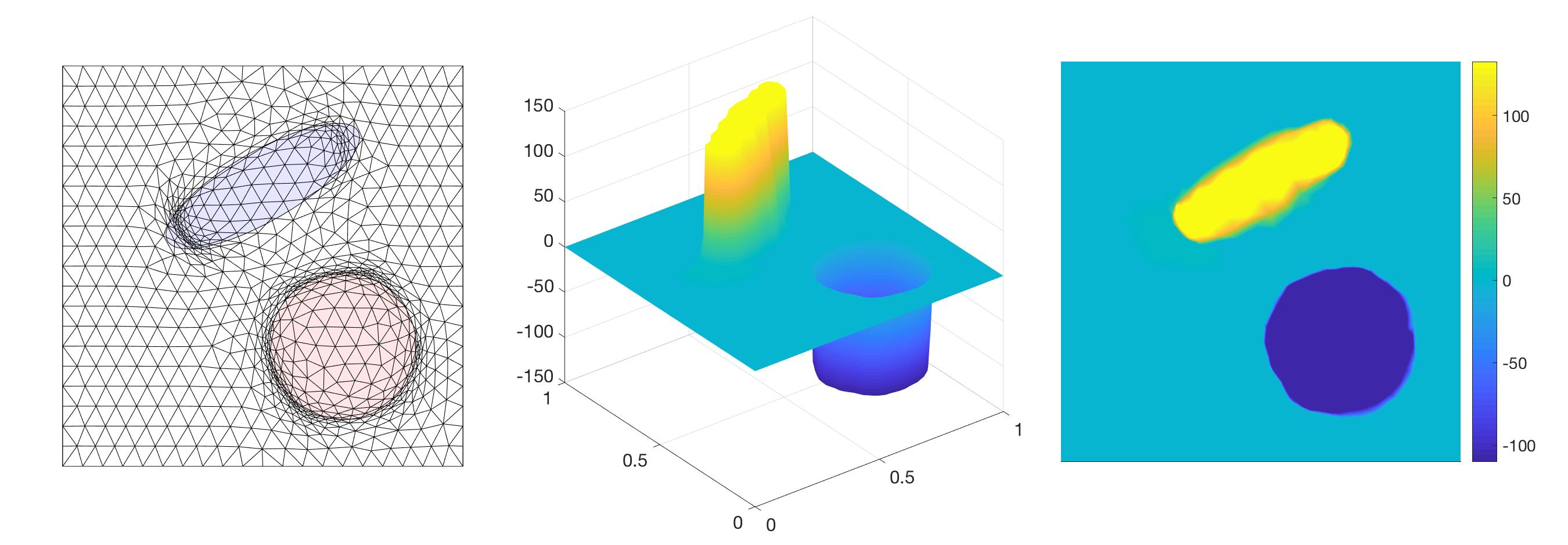}}
\caption{\label{fig:low noise Darcy} Four iterations of the hybrid IAS-remeshing algorithm. }
\end{figure} 

The results, together with the true source used for generation of the data are shown in Figure~\ref{fig:low noise Darcy}. 
The results show similar characteristics as in the tomography problem: The first iteration tends to overdensify the mesh near the discontinuities, after which the algorithm learns quickly the positions of the jumps and the mesh size decreases. The mesh sizes together with the computing times are again given in tabular form, see Table~\ref{tab:mesh Darcy}.

\begin{table}
\centerline{
\begin{footnotesize}
\begin{tabular}{c|cc|}
 & $n_t$ & time $[{\rm s}]$  \\
\hline
I  & 700 & 20.3  \\
II &5362 & 585.2 \\
III &1796 & 61.0  \\
IV &1588 & 33.5 \\
 \end{tabular}
\end{footnotesize}
}
\caption{\label{tab:mesh Darcy} The number of elements in the meshes, and the total computing time of the full updating of the estimates and the meshes.  The recomputing of the forward model is not included in the times.}
\end{table}

\section{Discussion and outlook}

Computational inverse problems for estimating distributed parameters are known to be sensitive to modeling errors, including those arising from insufficiently fine discretization of the underlying continuum model. This work develops further the ideas, outlined in the earlier work \cite{CCPS}, that the discretization is part of the inverse problem and should be decided concomitantly with the actual solution of the unknown parameter. Modeling errors are insignificant only if the discretization is fine enough to reduce the error below the noise level that needs to be taken into consideration in the solution of the problem, however, without a guidance of the refinement of the discretization, this may render the forward problem computationally too complex for efficient solution of the problem. In this work we are considering the anisotropic meshing connected to an anisotropic metric that is used to define tessellations allowing a numerical representation of the unknown and the solution of the forward model with sufficient accuracy to yield a modeling error dominated by the exogenous noise. Numerical tests show that the anisotropy reduces the computation times by a factor of 4-8 compared to isotropic selective meshing, and compared to uniform isotropic mesh refining, the reduction would be of orders of magnitude. The examples in this work are in two dimensions.  Generalization of the algorithm to three dimensions is possible albeit not without challenges. Furthermore, generalizations to non-linear inverse problems are underway and will be the topic of a separate contribution.

\section*{Acknowledgements}

The authors wish to acknowledge the partial support by the NSF, grants DMS 1951446 for Daniela Calvetti and DMS 2204618 for Erkki Somersalo. The work of Alberto Bocchinfuso was partly supported by the Great Lakes Energy Institute.

\end{document}